\documentclass[conf]{new-aiaa}
\usepackage[utf8]{inputenc}

\usepackage{graphicx}
\usepackage{amsmath}
\usepackage[version=4]{mhchem}
\usepackage{siunitx}
\usepackage{longtable,tabularx}
\usepackage{dsfont}
\usepackage{calc}
\usepackage{color}
\usepackage[ruled,vlined]{algorithm2e}
\usepackage{verbatim}
\newcommand{\norm}[1]{\left\lVert#1\right\rVert}
\newcommand{\abs}[1]{\left|#1\right|}
\newcommand{\menge}[1]{\left\lbrace#1\right\rbrace}

\newcommand{\intD}{\;\mathrm{d}}
\newcommand{\rom}[1]{\uppercase\expandafter{\romannumeral #1\relax}}

\newtheorem{theorem}{Theorem}
\newtheorem{proof}{Proof}

\newcommand{\inner}[1]{\left< #1 \right>}
\newcommand{\red}[1]{\textcolor{red}{#1}}

\title{A structure-preserving surrogate model for the closure of the moment system of the Boltzmann equation using convex deep neural networks}

\author{
Steffen Schotthöfer\footnote{Ph.D Candidate, steffen.schotthoefer@kit.edu},
Tianbai Xiao\footnote{Postdoc, tianbai.xiao@kit.edu},
Martin Frank\footnote{Professor, martin.frank@kit.edu}
}
\affil{Computational Science and Mathematical Methods Division,\\ Department of Mathematics, KIT Karlsruhe, \\ Bldg 20.30, Englerstraße 2, 76131 Karlsruhe, Germany}
\author{
Cory D. Hauck\footnote{Professor, hauckc@ornl.gov}
}
\affil{ Computer Science and Mathematics Division,\\ Oak Ridge National Laboratory, and
Department of Mathematics (Joint Faculty),  University of Tennessee, \\ 1 Bethel Valley Road, Bldg. 4100, Oak Ridge, TN 37831 USA }

\begin{document}

\maketitle

\begin{abstract}
Direct simulation of physical processes on a kinetic level is prohibitively expensive in aerospace applications due to the extremely high dimension of the solution spaces. 
In this paper, we consider the moment system of the Boltzmann equation, which projects the kinetic physics onto the hydrodynamic scale.
The unclosed moment system can be solved in conjunction with the entropy closure strategy. 
Using an entropy closure provides structural benefits to the physical system of partial differential equations. 
Usually computing such closure of the system spends the majority of the total computational cost, since one needs to solve an ill-conditioned constrained optimization problem. 
Therefore, we build a neural network surrogate model to close the moment system, which preserves the structural properties of the system by design, but reduces the computational cost significantly.
Numerical experiments are conducted to illustrate the performance of the current method in comparison to the traditional closure.

\end{abstract}


\section{Introduction}
In rarefied gas dynamics, the continuum assumption is no longer applicable and one has to rely on a more general description of the physical system at hand, which is given by kinetic equations such as the Boltzmann equation~\cite{BoltzmannApplications}. 
These equations arise in a variety of other applications as well, such as neutron transport~\cite{NeutronTransport}, radiative transport~\cite{Chahine1987FoundationsOR} and semiconductors~\cite{Markowich1990SemiconductorE}. 
The Boltzmann equation is a high dimensional integro-differential equation, with phase space dependency on time, space and particle velocity. This high dimensionality of the phase space presents a severe computational challenge for massive numerical simulations. 

Several methods for phase space reduction have been proposed to solve the Boltzman equation, including the discrete ordinate/velocity methods ~\cite{NeutronTransport,Camminady2019RayEM,xiao2020velocity} and moment methods~\cite{AlldredgeFrankHauck,AlldredgeHauckTits,GarretHauck,KRISTOPHERGARRETT2015573, Levermore}.
Discrete ordinate methods evaluate the velocity space at specific points, which yields a system of equations only coupled by the integral scattering operator. While computationally efficient, these methods suffers from numerical artifacts, which are called ray effects~\cite{Camminady2019RayEM}.
Moment methods eliminate the dependency of the phase space on the velocity variable by computing the moment hierarchy of the Boltzmann equation. Due to the structure of the advection term, the resulting moment system is unclosed. One distinguishes moment methods according to the modelling of their closure. 
The classical $P_N$ closure uses a simple truncation that results in a system of linear hyperbolic equations. The main drawback of this method is its numerical artifacts, specifically large oscillations of the particle density, which may result in negative particle concentrations. This effect is particularly present in the streaming particle regime~\cite{Brunner2002FormsOA}. 

A moment closure, which preserves important physical and mathematical properties~\cite{Levermore1996MomentCH} of the Boltzmann equation, is constructed by solving a convex constrained optimization problem, which is based on the entropy minimization principle~\cite{AlldredgeHauckTits, Levermore}.The method, which is commonly referred to as $M_N$ closure, is accurate in the diffusive limit~\cite{GOUDON2013579} and unlike the $P_N$ closure, this method is also accurate in the streaming limit~\cite{DUBROCA1999915}. Although the $M_N$ closure is methodically superior to the $P_N$ closure, it is by far more expensive to compute. Garret et al.~\cite{GarretHauck} have demonstrated, that in a high performance implementation, more than $80$\% of the computational time of the whole solver is required for the solution of the entropy minimization problem. This motivates the development of a neural network surrogate model to accelerate the $M_N$ closure. 

Several machine learning inspired methods have been proposed recently. 
The authors of~\cite{huang2021machine} close the moment system by learning the gradient of the highest order moment. In~\cite{Han21983}, the authors pursue two strategies. First, they use a encoder-decoder network to generate generalized moments and then learn the moment closure of system with its dynamics in mind. Second, they learn directly the correction term to the Euler equations. 
In~\cite{huang2020learning}, Galilean invariant machine learning methods for partial differential equations are developed using the conservation dissipation formalism. 
Using convolutional networks, a closure for the one dimensional Euler-Poisson system was constructed in ~\cite{bois2020neural}.
In~\cite{xiao2020using}, a dense neural network was used to model the deviation from the Maxwellian in the collision term of the Boltzmann equation.
The authors of~\cite{Maulik2020neural} use neural networks to reproduce physical properties of known magnetized plasma closures. 
In~\cite{ma20machine}, fully connected, dense and discrete Fourier transform networks are used to learn the Hammett-Perkins Landau fluid closure. 
Physics informed neural networks were employed  to solve forward and inverse problems via the Boltzmann-BGK formulation to model flows in continuum and rarefied regimes in~\cite{lou2020physicsinformed}, to solve the radiative transfer equation ~\cite{mishra2021physics} and the phonon Boltzmann equation in~\cite{li2021physicsinformed}. In~\cite{porteous2021datadriven}, the authors propose a data driven surrogate model of the minimal entropy closure using convex splines and empirically convex neural networks. To ensure convexity at the training data points, the authors penalize a non symmetric positive definite Hessian of the network output.  

The goal of this publication is to construct a deep neural network surrogate to model the entropy closure of the moment system of the Boltzmann Equation. The neural network and maps a given moment to its minimal mathematical entropy. In contrast to the work proposed in~\cite{porteous2021datadriven}, the neural network is input convex by design using the methods of Amos et al.~\cite{Amos2017InputCN}. By this ansatz, the learned closure automatically inherits all structural properties of the entropy closure for the moment system. The derivative of the network with respect to the moments maps to the corresponding optimal Lagrange multipliers of the entropy minimization problem. We train the neural network on the output, the Lagrange multiplier and additionally on the reconstructed moments, whereas in~\cite{porteous2021datadriven}, the authors train on the output, the reconstructed moments and the Hessian of the network.



The remainder of this paper is structured as follows.  In Section~\ref{sec_kinEq}, we briefly present the Boltzmann equation and its properties. We give a review to the moment method and the minimal entropy closure. 
In Section~\ref{sec_neuralEntropyClosure}, we  present the input convex neural network approach. Then  we construct the neural network closure as a surrogate model to the minimal entropy closure.
We show that the neural entropy closure has an entropy entropy-flux pair, fulfills a local dissipation law, is hyperbolicity preserving and invariant in range. Furthermore, it fulfills the H-theorem and the moment system conserves mass. We discuss the generation of training data and propose two algorithms for data generation. Lastly we describe the training and inference of the neural network on dimension reduced input data.
Section~\ref{sec_kineticSolver} describes the kinetic scheme for solving the moment system and its interplay with the neural network surrogate model.
Section~\ref{sec_numResults} contains a selection of $1D$ and $2D$ numerical test cases, where we demonstrate our findings.

\section{Kinetic equations and the moment method}\label{sec_kinEq}
Classical kinetic theory is profoundly built upon the Boltzmann equation,
which describes the space-time evolution of the one-particle kinetic density function $f(t, \mathbf{x},\mathbf{v})$ with $t>0$, $x\in \mathbf{X}\subset\mathbb R^3$, $\mathbf{v}\in\mathbf{V}\subset\mathbb R^3$ in a many-particle system,
\begin{equation}  \label{eq_boltzmann}
    \partial_t f+\mathbf{v} \cdot \nabla_{\mathbf{x}} f = Q(f).
\end{equation}
The left-hand side of the equation describes particle transport, while the right-hand side models collisions. If the particles only collide with a background material one can model this behavior with the  linear Boltzmann collision operator
\begin{equation}
    Q(f)(\mathbf{v})=\int_{\mathbf{V}} \mathcal B(\mathbf v_*, \mathbf v) \left[ f(\mathbf v_*)-f(\mathbf v)\right] d\mathbf v_*,
\end{equation}
where the collision kernel $\mathcal B(\mathbf v_*, \mathbf v)$ models the strength of collisions at different velocities. If the interactions among particles are considered, the collision operator becomes nonlinear. For example, the two-body collision results in
\begin{equation}
\label{eq_coll}
    Q(f)=\int_{\mathbf{V}} \int_{\mathcal S^2} \mathcal B(\cos \beta, |\mathbf{v}-\mathbf{v_*}|) \left[ f(\mathbf v')f(\mathbf v_*')-f(\mathbf v)f(\mathbf v_*)\right] d\mathbf \Omega d\mathbf v_*,
\end{equation}
where $\{\mathbf{v},\mathbf{v_*}\}$ are the pre-collision velocities of two colliding particles, and $\{\mathbf{v}',\mathbf{v_*}'\}$ are the corresponding post-collision velocities and $\mathcal S^2$ is the unit sphere.
The right-hand side is a fivefold integral, where $\beta$ is the so-called  deflection angle.

In the following, we use the notation
\begin{equation}
\inner{\cdot} = \int_{\mathbf{v}}\cdot \intD \mathbf{v}
\end{equation}
to define integrals over velocity space.

The Boltzmann equation is a first-principles model based on direct modeling.
It possesses some key structural properties, which are intricately related to the physical processes and its mathematical existence and uniqueness theory. We briefly review some of these properties, where we follow~\cite{AlldredgeFrankHauck,Levermore1996MomentCH}.
First, the time evolution of the solution is invariant in range, i.e. if $f(0,\mathbf{x},\mathbf{v})\in B\subset[0,\infty)$, then $f(t,\mathbf{x},\mathbf{v})\in B\subset[0,\infty)$ for all $t>0$. Particularly this implies non-negativity of $f$.
Second, if $\phi$ is a collision invariant fulfilling
\begin{align}\label{eq_collision_invariant}
        \inner{\phi Q(g)} = 0,\quad \forall g\in\mathrm{Dom}(Q),
\end{align}
the equation
\begin{align}
        \partial_t\inner{\phi f} + \nabla_{\mathbf x}\cdot\inner{\mathbf v\phi f} = 0
\end{align}
is a local conservation law. Third, for each fixed direction $\mathbf{v}$, the advection operator, i.e. the left-hand side term of Eq.~\eqref{eq_boltzmann}, is hyperbolic in space and time. Forth, let $D\subset\mathbb{R}$. There is a twice continuously differentiable, strictly convex function $\eta: D\rightarrow\mathbb{R}$, which is called kinetic entropy density.  It has the property
\begin{align}
        \inner{\eta'(g)Q(g)} \leq 0, \quad \forall g\in\mathrm{Dom}(Q) \ \mathrm{s.t.} \ \mathrm{Im}(g)\subset D.
\end{align}
Applied to Eq.~\eqref{eq_boltzmann}, we get the local entropy dissipation law
\begin{align}
    \partial_t\inner{\eta(f)} + \nabla_{\mathbf x}\cdot\inner{\mathbf v\eta(f)} \leq 0.
\end{align}
Usually we set $D=B$. Lastly, the solution $f$ fulfills the H-theorem.

The Boltzmann equation is an integro-differential equation model defined on a seven-dimensional phase space. With the nonlinear five-fold integral, it is  challenging to solve accurately and efficiently. 
The well known moment methods encode the velocity dependence of the Boltzmann equation in a moment vector $u\in\mathbb{R}^{N+1}$ of order $N$~\cite{AlldredgeFrankHauck, Levermore1996MomentCH}.
The moments are calculated with respect to a vector of basis functions $m(\mathbf{v})\in\mathbb{R}^{N+1}$ by integrating over the velocity space,
\begin{align}\label{eq_momentDef}
 u(t,\mathbf x)=\inner{m f}.
 \end{align} 
Common choices for the basis functions are monomials and spherical harmonics, depending on the application. Typically, they include the collision invariants defined in Eq.~\eqref{eq_collision_invariant}.
The moment vector satisfies the system of transport equations
\begin{align}\label{eq_momentRT}
\partial_t u(t,\mathbf{x}) + \nabla_\mathbf{x}\cdot\inner{\mathbf{v} m(\mathbf{v})f}= \inner{m(\mathbf{v})Q(f)},
\end{align}
which is called moment system.
By construction, the advection operator depends on $f$ and thus, the moment system is unclosed. Moment methods aim to find a meaningful closure for this system. 
Since the kinetic equation dissipates entropy and fulfill a local entropy dissipation law, one can close the system by choosing the reconstructed kinetic density $\tilde{f}$ out of all possible functions $g$ that fulfill $u(t,\mathbf{x})=\inner{mg}$ as the one with minimal entropy $h$.  The minimal entropy closure can be formulated as a constrained optimization problem for a given set of moments $u$.
\begin{align}\label{eq_entropyOCP}
\begin{aligned}
 h(u) = \inner{\eta(f_u)} &= \min_{g\in F} \inner{\eta(g)} \\
 \text{ s.t. }  u&=\inner{m g}
\end{aligned}
\end{align}
where $F=\menge{g : \text{Range}(g)\subset D}$.
The set $\mathcal{R}$ of all moment vectors $u$ for which a solution exists is called the set of realizable moments. These are all moments that belong to a kinetic density function $f$ that fulfills the invariant range property of the kinetic equation. This condition is enforced by the constraint of the optimization problem of Eq.~\eqref{eq_entropyOCP}, which is called realizability constraint.
For all moment vectors $u$ for which a solution exists, this solution is unique and of the form 
\begin{align}\label{eq_entropyRecosntruction}
    {f}_u = \eta'_*(\alpha(u)\cdot m),
\end{align}
where the Lagrange multiplier $ \alpha:\mathbb{R}^{N+1}\rightarrow\mathbb{R}^{N+1} $ maps $u$ to the solution of the dual problem
\begin{align}\label{eq_entropyDualOCP}
    \alpha_u(u) =  \underset{\alpha\in\mathbb{R}^{N+1}}{\text{argmax}} \menge{ \alpha\cdot u - \inner{\eta_*(\alpha\cdot m)}}
\end{align}
and $\eta_*$ is the Legendre dual of $\eta$.
Inserting ${f}_u$ into the moment system \eqref{eq_momentRT} gives a closed system of equations in space and time. The function
\begin{align}\label{eq_entropyFunctionalH}
h(u) =  \alpha_u\cdot u - \inner{\eta_*(  \alpha_u\cdot m)}
\end{align}
is  convex in $u$ and a suitable entropy functional for the moment system in Eq.~\eqref{eq_momentRT}, see~\cite{AlldredgeFrankHauck} .
Furthermore,  the derivative of $h$ recovers the optimal Lagrange multipliers of Eq.~\eqref{eq_entropyDualOCP},
\begin{align}\label{eq_derivH}
\frac{\intD}{\intD u}h =  \alpha_u
\end{align}
This minimal entropy closure  also conserves the above listed structural properties of the Boltzmann equation . We present the above properties for the moment system for the sake of completeness, where we follow~\cite{AlldredgeFrankHauck,Levermore}. First, the invariant range property of the solution $f$ translates to the set of realizable moments $\mathcal{R}$, which can be described as the set of moments corresponding to a kinetic density with range in $B$. One demands that $u(t,\mathbf{x},\mathbf{v})\in\mathcal{R}$ for all $t>0$. Second, if a moment basis function $m_i(\mathbf{v})$ is a collision invariant, then 
 \begin{align}\label{eq_momentRT_lhS}
        \partial_t u(t,x) + \nabla_{\mathbf{x}}\cdot\inner{\mathbf{v} m f_u}=0,
\end{align}
is a local conservation law. Third, one can write Eq.~\eqref{eq_momentRT} as a symmetric hyperbolic conservation law in $\alpha_u$. Forth,  for  $u\in\mathcal{R}$, $h(u)$ and $j(u)=\inner{\mathbf{v}\eta(f_u)}$ is a suitable entropy and entropy-flux pair compatible with the advection operator $\inner{\mathbf{v} mf_u}$ and yield a semi-discrete version of the entropy dissipation law. 
\begin{align}
    \partial_th(u) +\nabla_{\mathbf{x}} j(u)=h'(u)\cdot\inner{mQ(f_u(\alpha_u)}\leq 0
\end{align}
Note that convexity of $h(u)$ is crucial for the entropy dissipation property. Lastly, the moment system fulfills the H-theorem. 

A numerical method to solve the moment system therefore consists of an iterative discretization scheme for the moment system~\eqref{eq_momentRT} and a   Newton optimizer for the dual minimal entropy optimization problem in Eq.~\eqref{eq_entropyDualOCP}. The former scheme can be a finite volume or discontinuous Garlerkin scheme, for example. 
The drawback of the method is the high computational cost associated with the Newton solver. The optimization problem in Eq.~\eqref{eq_entropyDualOCP} needs to be solved in every grid cell at every time step of the PDE solver. The computational effort to solve the minimal entropy optimization problem grows over-proportionate with the dimension $N$ of the moment basis $m$. Using three basis functions, the optimizer requires $80$\% of the computation time and $87$\% when using seven basis functions, as Garrett et al. have demonstrated in a computational study~\citep{KRISTOPHERGARRETT2015573}. 
Furthermore, the optimization problem is ill-conditioned, if the moments $u$ are near the boundary of the realizable set $\mathcal{R}$~\cite{AlldredgeHauckTits}. At the boundary $\partial\mathcal{R}$, the Hessian of the objective function becomes singular and the kinetic density $f_u$ is a sum of delta functions~\cite{Curto_recursiveness}.

\section{Neural entropy closures}\label{sec_neuralEntropyClosure}
In the following, we propose a neural network to solve the minimal entropy closure problem of Eq.~\eqref{eq_entropyDualOCP}, which preserves the intricate structure of the moment system and the Boltzmann equation. The advantage of a neural network compared to a Newton optimizer is a significantly lower computational expense. Instead of inverting a possibly near singular Hessian multiple times,  a trained neural network only needs a comparatively small amount of fast tensor operations to compute the closure.

 \subsection{Design of the neural network closure}


A neural network $\mathcal{N}_\theta: \mathbb{R}^n\mapsto\mathbb{R}^m$ is a parameterized mapping from an input $x$ to the network prediction $\mathcal{N}_\theta(x)$. In this work, we focus on dense neural networks. 
A multi-layer neural network $\mathcal{N}_\theta$ is a concatenation of non-linear (activation) functions $f_k$ applied to weighted sums of the previous layer's output $z_{k-1}$. An $M$ layer network can be described in tensor formulation as follows.
\begin{align}\label{eq_nnKons}
z_k &= f_k(W_kz_{k-1} +b_k), \qquad	k = 1,\dots,M \\
x &= z_0,  \\
\mathcal{N}_\theta(x) &= z_M
\end{align}
where $W_k$ is the weight matrix of layer $k$ and $b_k$ the corresponding bias vector. In the following, we denote the set of all trainable parameters of the network, i.e. weights and biases by $\theta$.
With the correct trainable parameters, a neural network is capable of approximating any given target function with arbitrary accuracy~\cite{cybenko}, however it is not trivial how to choose them. 
Usually, one chooses a set of training data points  $T=\menge{(x_j,y_j)}_{j\in J}$ to input in a loss function, for example the mean squared error
\begin{align}
L(x,y;\theta) = \frac{1}{\vert J \vert}\sum_{j\in J} \norm{y_j -\mathcal{N}_\theta(x_j)}^2_2.
\end{align} 
Then one can set up the process of finding suitable weights, called training of the network, as an optimization problem
\begin{align}
\min_\theta \,  &L(x,y;\theta)
\end{align}
The optimization is often carried out with gradient-based algorithms, such as stochastic gradient descent~\cite{robbins1951} or related methods as ADAM~\cite{kingma2017adam}, which we use in this work.

The proposed neural network based closure for the minimal entropy problem, called neural entropy closure, is a network that maps a realizable moment vector $u\in\mathcal{R}$ to the  entropy functional $h(u)$  in Eq.~\eqref{eq_entropyFunctionalH}, which is the  objective functional of the entropy optimization problem in Eq.~\eqref{eq_entropyOCP}. 
Assuming the neural network is trained, we have the following relations,
\begin{align}
    h_\theta :=& \mathcal{N}_\theta \approx h(u),\\
    \alpha_\theta :=& \frac{\intD}{\intD u}\mathcal{N}_\theta \approx \frac{\intD}{\intD u}h = \alpha_u,\\
    f_\theta :=&\eta'_*(\frac{\intD}{\intD u}\mathcal{N}_\theta\cdot m) \approx  \eta'_*(\alpha_u\cdot m) = f_u, \\
    u_\theta :=& \inner{m\eta'_*(\frac{\intD}{\intD u}\mathcal{N}_\theta\cdot m)} \approx  \inner{m\eta'_*(\alpha_u\cdot m)} = u, \label{eq_u_theta}
\end{align}
by using Eq.~\eqref{eq_entropyRecosntruction}, Eq.~\eqref{eq_derivH} and the definition of the moment.
These relations give a set of possible loss functions for the network training and a workflow for the neural network accelerated kinetic solver. We choose $h$, $\alpha_u$ and $u$ as network outputs, since $f_u$ is dependent on the velocity variable $\mathbf{v}$. The loss function of the network then becomes
\begin{align}\label{eq_nwLoss}
  L(u,\alpha_u,h;\theta) =  \frac{1}{\vert J\vert}\sum_{j\in J} \norm{h(u_j) -\mathcal{N}_\theta(u_j)}^2_2+\norm{\alpha_u(u_j)- \alpha_\theta(u_j)}^2_2 +\norm{u_j- u_\theta(u_j)}^2_2.
\end{align}
We demand that the neural entropy closure preserves the structure of the moment system of the Boltzmann equation, which is described in Section~\ref{sec_kinEq}. 
The invariant range property of $f_\theta$ depends solely on the range of $\eta_*'$. Popular choices for the kinetic entropy density $\eta$ are Maxwell-Boltzmann, Bose Einstein or Fermi-Dirac. In this work, we focus on the Maxwell-Boltzmann entropy, which has the following definition, Legendre dual and derivative.
\begin{align}
    \eta(z) &= z\ln(z)-z, \qquad z\in D=\mathbb{R}_+ \\
    \eta'(z) &= \ln(z), \qquad\qquad\  z\in D=\mathbb{R}_+ \\
    \eta_*(y) &= \exp(y),  \qquad\quad\ \ y\in\mathbb{R}\\
    \eta_*'(y) &= \exp(y),  \qquad\quad\ \ y\in\mathbb{R}
\end{align}
By construction, the neural entropy closure is of invariant range, since $f_\theta(v)=\exp(\frac{\intD}{\intD u}\mathcal{N}_\theta\cdot m(v))>0$.
Interchanging the entropy functional by a neural network does not affect the conservation property of the moment system.
Consider the hyperbolicity requirement. In order to define the Legendre dual of $h$, it must be convex. In the proof of the hyperbolicity property, which is conducted in \cite{Levermore1996MomentCH} for $\alpha_u$ and $u$ as the system variable, $h''$ (respectively $h_*''$) must be symmetric positive definite. As a consequence, $h(u)$ and therefore the neural network $N_\theta(u)$ must be strictly convex in $u$. Strict convexity of the entropy functional $h$ is the crucial requirement for the related properties entropy dissipation and the H-theorem~\cite{Levermore1996MomentCH}.

Convex neural networks have been inspected in \cite{Amos2017InputCN}, where the authors propose several deep neural networks that are strictly convex with respect to the input by design. The design is led by the following principles~\cite{boyd_vandenberghe_2004}. First, a positive sum of convex functions is convex. Second, let $f:\mathbb{R}^n\rightarrow\mathbb{R}$ be the concatenation of the functions $h:\mathbb{R}^k\rightarrow\mathbb{R}$ and $g:\mathbb{R}^n\rightarrow\mathbb{R}^k$. Then $f(x) = h(g(x))$ is convex, if $h$ is convex, $h$ is non-decreasing in each argument and all $g_{i=1,\dots,k}$ are convex.
Applying these conditions to the definition of a layer of a neural network, Eq.~\eqref{eq_nnKons}, we get that all entries of the weight matrix $W_k$ must be  positive in all layers except the first. Furthermore, the activation function of each layer must be strictly convex. However, in practice it turns out that very deep networks with positive weights have difficulties to train. The authors therefore modifies the definition of a hidden layer in Eq.~\eqref{eq_nnKons} to
\begin{align}\label{eq_ICnnKons}
z_{k} &= f_k(W_k^z z_{k-1} + W_k^x x + b_k^z), \qquad k=2,\dots, M, \\
z_{k} &= f_k(W_k^x x + b_k^z), \qquad\qquad\qquad\, k=1,
\end{align}
where $W_k^z$ must be non-negative, and $W_k^x$ may attain arbitrary values. This approach is related to the idea of residual neural networks~\cite{He2016DeepRL}. We choose the strictly convex softplus function $ {\displaystyle \ln \left(1+e^{x}\right)}$ as layer activation function for $f_k, k=1,\dots,M-1$
and a linear activation for the last layer, since we are dealing with a regression task.

\subsection{Training Data}\label{subsec_Data}
The inputs to the neural entropy closure are realizable moments $u\in\mathcal{R}$ and it is trained on the corresponding Lagrange multipliers $\alpha\in\mathbb{R}^{N+1}$ and the values of the entropy functional $h(u,\alpha)$ defined in Eq.~\eqref{eq_entropyDualOCP}. The last output of the network is u itself, where is measured, how well the network can reconstruct the input moment, see Eq.~\eqref{eq_u_theta}.
The training datset therefore consists of triplets $(h(u,\alpha_u,\alpha_u,u)$, which have the relations described in Eq.~\eqref{eq_entropyDualOCP} and Eq.~\eqref{eq_entropyFunctionalH}.
Unfortunately, $\mathcal{R}$ is unbounded, since $u_0 = \inner{f} \geq 0$, if $f\geq 0$ and thus all non-negative moments of order zero are realizable. A strategy to characterize the realizable set, is to consider normalized moments \cite{Kershaw1976FluxLN}
\begin{align}
    \overline{u} &= \frac{u}{u_0} = [1,\frac{u_1}{u_0},\dots,\frac{u_N}{u_0}]^T\in\mathbb{R}^{N+1},\\
    \overline{u}^r &= [\overline{u}_1^r,\dots,\overline{u}_N^r]^T\in\mathbb{R}^N,
\end{align}
where we call $\overline{u}$ the normalized moments and  $\overline{u}_r$ the reduced normalized moments. Analogously, we define the set of normalized realizable moments $\overline{\mathcal{R}}$ and the reduced normalized realizable moments as  $\overline{\mathcal{R}}^r$
\begin{align}
    \overline{\mathcal{R}} &=  \menge{u\in\mathcal{R}: u_0 =1}\subset\mathbb{R}^{N+1},\\
    \overline{\mathcal{R}}^r &=  \menge{\overline{u}^r\in\mathbb{R}^{N} [1,\overline{u}^T]^T\in\overline{\mathcal{R}}}\subset\mathbb{R}^{N}.
\end{align}
Both, $\overline{\mathcal{R}}$ and $\overline{\mathcal{R}}^r$ are bounded and convex~\cite{Kershaw1976FluxLN,Monreal_210538}. 
In \cite{Kershaw1976FluxLN}, the authors have derived expressions for  $\overline{\mathcal{R}}^r$ up to moment order $N=4$ and one spatial dimension, i.e. $V,X\subset\mathbb{R}^1$ for monomial basis functions.
\begin{subequations}\label{eq_kershaw}
\begin{align}
    1&\geq \overline{u}^r_1 \geq -1 \\
    1&\geq \overline{u}^r_2 \geq (\overline{u}^r_1)^2 \\
    \overline{u}^r_2 -\frac{(\overline{u}^r_1-\overline{u}^r_2)^2}{1-\overline{u}^r_1} &\geq \overline{u}^r_3 \geq - \overline{u}^r_2 + \frac{(\overline{u}^r_1+\overline{u}^r_2)^2}{1+\overline{u}^r_1} \\
    \frac{(\overline{u}^r_2)^3-(\overline{u}^r_3)^2+2 \overline{u}^r_1 \overline{u}^r_2 \overline{u}^r_3}{\overline{u}^r_2-(\overline{u}^r_1)^2}&\geq \overline{u}^r_4 \geq  \overline{u}^r_2 - \frac{(\overline{u}^r_1 -\overline{u}^r_3)^2}{(1-\overline{u}^r_2)}
\end{align}
\end{subequations}
where equality indicates the boundary $\partial\overline{\mathcal{R}}^r$. Here, the basis functions are monomials.
The realizable set for higher order moments can be characterized using the more general results in~\cite{Curto_recursiveness}.
This characterization enables a sampling strategy for training data, where one first samples uniformly in $\overline{\mathcal{R}}^r$ and afterwards solves the dual entropy minimization problem in Eq.~\eqref{eq_entropyDualOCP} using $\overline{u}$, to compute $\alpha_u$ and $h(u,\alpha_u)$, see Algorithm~\ref{alg_samplingUniform}.
\begin{figure}
\centering
    \begin{minipage}{0.49\textwidth}
      \begin{algorithm}[H]\label{alg_samplingUniform}
\SetAlgoLined
\KwResult{$\menge{(\overline{u}_j,\alpha_j,h_j)}_{j\in J}$}
 Set $N$, $\vert J\vert$ \;
 Set the training accuracy $\tau>0$\;
 Set the boundary distance $\delta>0$\;
 \For {$j = 1,\dots,\vert J\vert$}{
 Sample $\overline{u}_j$ using Eq.~\eqref{eq_kershaw} such that dist$(u_j,\partial\overline{\mathcal{R}})>\delta$\;
 Solve Eq.~\eqref{eq_entropyDualOCP} for $\alpha_{u,j}$ using a Newton optimizer with tolerance $\tau$\;
 Compute $h_j$ using Eq.~\eqref{eq_entropyFunctionalH}\;
 }
 \caption{Data generation with uniform sampling of $\overline{\mathcal{R}}$ }
\end{algorithm}
   \end{minipage}
   \begin{minipage}{0.49\textwidth}
       \begin{algorithm}[H]\label{alg_samplingUniformAlpha}
\SetAlgoLined
\KwResult{$\menge{(\overline{u}_j,\alpha_j,h_j)}_{j\in J}$}
 Set $N$, $\vert J\vert$ \;
 Determine $A\subset\mathbb{R}^N$\;
 \For {$j = 1,\dots,\vert J\vert$}{
 Sample $\alpha^r_j$ such that $\alpha\in A$\;
 Compute $\alpha_{0,j}$ using Eq.~\eqref{eq_alpha0_recons}\;
 Set $\alpha_j = [\alpha_{0,j},{\alpha^r_j}^T]$\;
 Compute $\overline{u_j}$ using Eq.~\eqref{eq_entropyRecosntruction}\;
 Compute $h_j$ using Eq.~\eqref{eq_entropyFunctionalH}\;
 }
 \caption{Data generation with uniform sampling of Lagrange multipliers}
\end{algorithm}
   \end{minipage}
\end{figure}
The algorithm has some drawbacks. First, we do not solve Eq.~\eqref{eq_entropyDualOCP} exactly, but numerically with a tolerance $\tau$. As a result, the neural network accuracy is bounded from below by $\tau$. In practise however, a Newton optimizer is more accurate than the training accuracy of a neural network.  Second, 
the optimization problem Eq.~\eqref{eq_entropyDualOCP}  becomes ill-conditioned near the boundary $\partial\overline{\mathcal{R}}^r$ and singular at $\partial\overline{\mathcal{R}}^r$, so we can not generate training data at the boundary using this approach. Lastly, necessary and sufficient conditions for  $\partial\overline{\mathcal{R}}^r$ are only derived for $\mathbf{V},\mathbf{X}\subset\mathbb{R}^1$. Necessary conditions for $\mathbf{V},\mathbf{X}\subset\mathbb{R}^d$, $d=2,3$ are stated in~\cite{Monreal_210538} for $N\leq 2$. In $1$D $M_1$, they simplify to 
\begin{align}
    \norm{\overline{u}^r}\leq 1.
\end{align}
A different strategy for generating training data is to sample the Lagrange multipliers $\alpha_u$ and then reconstruct $u$ with Eq.~\eqref{eq_entropyRecosntruction} and Eq.~\eqref{eq_momentDef}. 
In order to generate normalized moments, we have to construct the image of $\overline{\mathcal{R}}$ which is a $N$ dimensional hyperplane. For a monomial basis and the Maxwell-Boltzmann entropy, we get the dependent Lagrange multiplier by the definition of the moment of order zero.
\begin{align}
1 = \inner{m_0 \eta_*'(\alpha^Tm)} = \inner{\exp(\alpha\cdot m)} =\inner{\exp(\alpha^r\cdot m^r)\exp(\alpha_0)},
\end{align}
which we can transform to 
\begin{align}\label{eq_alpha0_recons}
    \alpha_0 = -\ln(\inner{\exp(\alpha^r\cdot m^r)}),
\end{align}
where $\alpha^r = [\alpha_1,\dots,\alpha_N]^T$ and $m^r = [m_1,\dots,m_N]^T$.
For numerical stability, one needs to restrict the sampled $\alpha^r$ to a bounded set $A\subset\mathbb{R}^N$, since  as $u$ approaches $\partial\mathcal{R}$, $\norm{\alpha}\rightarrow\infty$. Algorithm~\ref{alg_samplingUniformAlpha}  summarizes this approach. 
Computational resources for both sampling algorithms can be found in the open source framework KiT-RT, see~\cite{KITRT}.

\subsection{Inference and training of the neural entropy closure}
In Section~\ref{subsec_Data}, we have discussed a normalization and a dimension reduction of the input data, which enables the neural network to train and execute on a bounded and convex set of possible inputs and therefore drastically increases its reliability. Each unseen data point is in the convex hull of seen data points. Normalization of $\mathcal{R}$ allows for a dimension reduction, thus the neural entropy closure takes normalized reduced moments $\overline{u}^r$ as . One can think of the constant zero order moment being encrypted in the bias term of the first layer. The primary output $\mathcal{N}_\theta(\overline{u}^r)$ is trained on $h(\overline{u},\alpha)$, but the input derivative $\partial_{\overline{u}^r}$ now only reveals $\alpha_\theta^r$. The missing value of $\alpha_0$ can be computed using Eq.~\eqref{eq_alpha0_recons}. The third network output, $\overline{u}_\theta$ is computed using Eq.~\eqref{eq_alpha0_recons}, Eq.~\eqref{eq_entropyRecosntruction} and Eq.~\eqref{eq_momentDef}. In order to train the network on $\alpha$, respectively $\alpha^r$, we compute the sensitivity of the input derivative of the neural network, which can be thought of training in Sobolev norm. This concept has been investigated in~\cite{SobolevTraining}, where the authors find that training in Sobolev norm increases training accuracy and data efficiency. 

After training, the neural network has the task to compute the minimal entropy closure for $u\in\mathcal{R}$.
Inference of a normalized moment $\overline{u}$ with $\mathcal{N}_\theta$ produces $\alpha_\theta(\overline{u})$. 
By considering again the definition of moment zero
\begin{align}
     1 &= \inner{m \exp(\alpha_\theta(\overline{u}))}  \\
    \iff   u_0  &= \inner{m \exp(\ln(u_0)\exp(\alpha_\theta(\overline{u}))}
\end{align}
we can deduce for $u_0>0$
\begin{align}
    \alpha_\theta(u) = [\alpha_\theta(\overline{u})_0 + \ln(u_0),\alpha_\theta(\overline{u})_1,\dots,\alpha_\theta(\overline{u})_N]^T.
\end{align}
Then $\alpha_\theta(u)$ minimizes the entropy functional for $u$, if and only if $\alpha_\theta(\overline{u})$ minimizes the entropy functional for $\overline{u}$. \\
Note, that solving the normalized minimal entropy closure, respectively inference of the network with the normalized moments can be interpreted as replacing the kinetic entropy density $\eta$ by the normalized kinetic entropy density 
\begin{align}
    \overline{\eta}(f) = \eta\left(\frac{f}{\inner{f}}\right)
\end{align}
which is still a strictly convex function and therefore fulfills the hyperbolicity, entropy-dissipation and H-theorem properties.

\section{Kinetic Scheme for the Moment Method}\label{sec_kineticSolver}
In this section, we briefly describe the kinetic scheme that is used to solve the moment method described in Section~\ref{sec_kinEq}. 
 The idea of a kinetic scheme~\cite{GarretHauck,AlldredgeHauckTits,KRISTOPHERGARRETT2015573} is to combine a numerical scheme for hyperbolic conservation laws with an according velocity space discretization.
 
 A Quadrature approximation of the velocity integral is needed in several parts of the kinetic scheme, namely the computation of the numerical flux and the collision operator $Q$.
\begin{align}
\inner{f(\mathbf{v})} \approx \sum_{q=1}^{n_q} w_q f(\mathbf{v_q}),
\end{align}
where $w_q$ are the quadrature weights and $v_q$ the corresponding quadrature points.
Depending on the dimension $d$ of $\mathbf{V}$, we use different parametrizations and quadrature rules~\cite{atkinson1982,atkinson2012spherical}, see Table~\ref{tab_quad}.
\begin{table}[htbp]
	\caption{Parametrization of $\mathbf{V}$ and used quadrature rule} 
	\centering
	\begin{tabular}{|l|l|l|l|l|l|l|} 
		\hline
		$d$ & $\mathbf{V}$& parametrization of $\mathbf{v}$ & Quadrature Rule \\ 
		\hline 
		$1$ & $[-1,1]$ & $\mathbf{v}= \mu$ & Legendre \\ 
		\hline 
		$2$ & $[-1,1]\times[0,2\pi)$ & $\mathbf{v}= \left( \begin{array}{c} \sqrt{1-\mu}\cos(\phi) \\ \sqrt{1-\mu}\sin(\phi) \end{array}\right)$ & Tensorized Gauss-Legendre, projected \\ 
		\hline 
		$2$ & $[-1,1]\times[0,2\pi)$ & $\mathbf{v}= \left( \begin{array}{c} \sqrt{1-\mu}\cos(\phi) \\ \sqrt{1-\mu}\sin(\phi)\\ \mu \end{array}\right)$ & Tensorized Gauss-Legendre \\ 
		\hline
	\end{tabular} 
	\label{tab_quad}
\end{table}

The spatial discretization is implemented with a finite volume method using a quadrilateral grid in the $2$ dimensional test cases. The notification is similar for the $1$ dimensional test case. 
We consider the flux function
\begin{align}\label{eq_flux}
   F(f,\mathbf{v}) =  \nabla_{\mathbf{x}}\cdot\inner{\mathbf{v}m(\mathbf{v})f_u}
\end{align}
We denote the moment vector averaged over control volume $i$ as $u^i$ and average Eq.~\eqref{eq_boltzmann} over the area $A_i$ of the control volume. The flux function $F(f,v)$ is approximated using an out of the box numerical flux $F_n(f_{u,i},f_{u,j},\mathbf{v})$ at the face of control volume $i$ and its neighbor $j$. This yields the kinetic flux $G(i,\mathbf{v})$ at control volume $i$
\begin{align}
   G(i,\mathbf{v}) = \frac{1}{A_i} \inner{\sum_{j\in\mathcal{N}(i)}F_{up}(i,j,\mathbf{v})l_{i,j}}, 
\end{align}
where $\mathcal{N}(i)$ are all neighboring grid cells of control volume $i$ and $l_{i,j}$ is the length of the face  of control volume $i$ and its neighbor $j$.
Finally, we approximate the integral with a suitable integration rule and get
\begin{align}
   G(i) = \frac{1}{A_i} \sum_{q = 1}^{n_q} w_q \sum_{j\in\mathcal{N}(i)}F_{up}(i,j,\mathbf{v}_q)l_{i,j}.
\end{align}
This means, that we have to evaluate the numerical flux for each quadrature point and for each equation in the moment system, which is quite costly but easy to parallelize.

We discretize the moments of the collision operator similarly. Note that the collision integral simplifies, if the chosen moments are collision invariants. The collision operator $Q(f_{u,i})$, respectively its moments $\inner{m(v)Q(f_{u,i})}$ at control volume $i$ can be computed using the above described quadrature rule
\begin{align}
     Q(f_{u,i}(\mathbf{v}))\approx Q_n(f_{u,i}(\mathbf{v})) = \sum_{q=1}^{n_q} w_q \mathcal B(\mathbf {\mathbf{v}_*}_q, \mathbf {v}) \left[ f_{u,i}(\mathbf{v}_{*q})-f_{u,i}(\mathbf{v})\right],
\end{align}
and
\begin{align}
   \frac{1}{A_i} \inner{m(\mathbf{v})Q(f_{u,i})} \approx  D(i) =\frac{1}{A_i} \sum_{q=1}^{n_q}w_q m(\mathbf{v}_q)Q_n(f_{u,i}(\mathbf{v}_q)),
\end{align}
which we denote by $D(i)$.
If there is a source term in addition to the collision term, it can be discretized analogously. \\
The temporal discretization can be done by any time-stepping scheme for hyperbolic partial differential equation. 
We consider the spatial discretization of the kinetic equation as a system of ordinary differential equations
\begin{align}
   \partial_t u^i =  - G(i;t) + D(i;t),
\end{align}
which is then discretized by an explicit Euler method with step size $\Delta t$
\begin{align}
    u^{t+\Delta t,i} = u^{t,i} + \Delta t\left(  D(i;t) - G(i;t)\right).
\end{align}
Recall, that we need to compute $f_u$ using either a numerical Newton solver for the optimization problem in Eq~\eqref{eq_entropyDualOCP} or by evaluating the neural entropy closure $N_\theta$. Either way, we introduce an error and thus the reconstructed kinetic density $f_{u^t}$ does not match exactly the moments $u^t$.  As a consequence, the updated moment of the next time step $u^{t+\Delta t}$ is not necessarily realizable~\cite{realizability_Kusch_Frank, AlldredgeHauckTits, KRISTOPHERGARRETT2015573}. In~\cite{AlldredgeHauckTits}, the authors propose a time-step restriction to preserve realizability of the finite volume updatem, whereas in~\cite{realizability_Kusch_Frank}, the authors change the current moment $u^{t,i}$ to $\inner{m f_u}$. This procedure yields an exact pair of moments and Lagrange multipliers, but comes to the cost of losing the conservation property of the scheme. Computational resources for the kinetic solver can be found in the open source kinetic solver Kinetic.jl~\cite{xiao2021julia,KitML}.
The neural networks are implemented in tensorflow 2.0~\cite{tensorflow2015-whitepaper} and the computational resources for the networks in this paper can be found in the repository~\cite{NeuralEntropy}.

\section{Numerical Results} \label{sec_numResults}
In this section, we present numerical results  and  investigate the performance of the neural entropy closure. First, we compare the sampling strategies in  Algorithm~\ref{alg_samplingUniform} and Algorithm~\ref{alg_samplingUniformAlpha}. Then we train the networks using the generated data. Finally, we employ the network in the kinetic solver and compare the results with the benchmark solution, where the minimal entropy problem is solved with a Newton solver.

\subsection{Comparison of sampling strategies}
\begin{figure}
    \centering
   \begin{minipage}{0.49\textwidth}
   \centering
       \includegraphics[width=\textwidth]{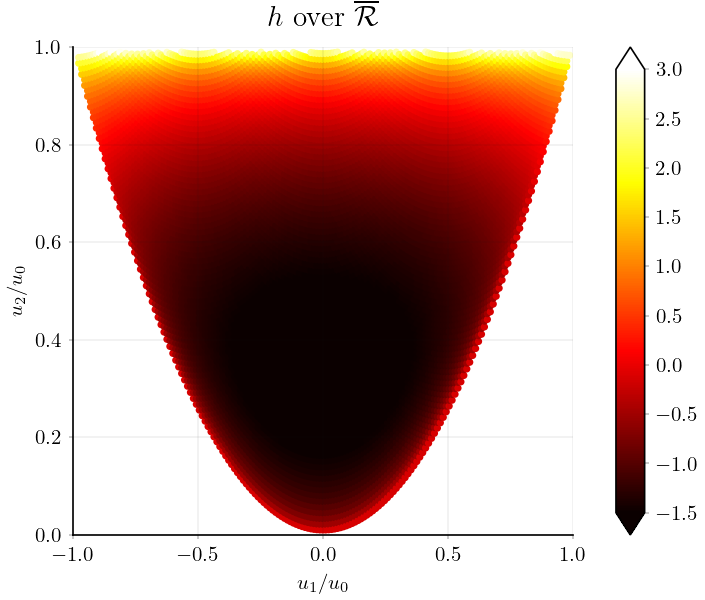}
               {\\a) Sampling with Algorithm~\ref{alg_samplingUniform}}
   \end{minipage}
   \begin{minipage}{0.49\textwidth}
     \centering
        \includegraphics[width=\textwidth]{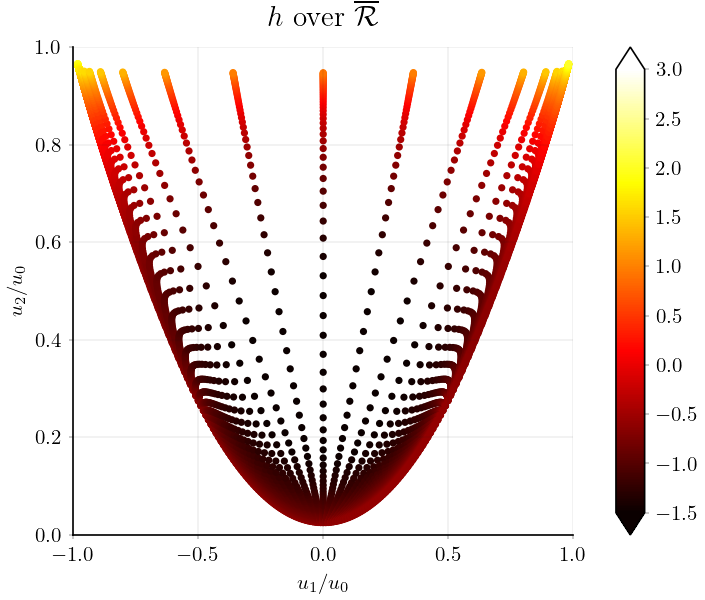}
        {\\b) Sampling with Algorithm~\ref{alg_samplingUniformAlpha}}
   \end{minipage}
    \caption{Scatter plots of $\overline{R}^r$ for  $N=2$ with data generated from different strategies. Color indicates the value of the minimum point of the entropy functional.}
    \label{fig_dataSampling}
\end{figure}
\begin{figure}
    \centering
    \begin{minipage}{0.49\textwidth}
    \centering
        \includegraphics[width=\textwidth]{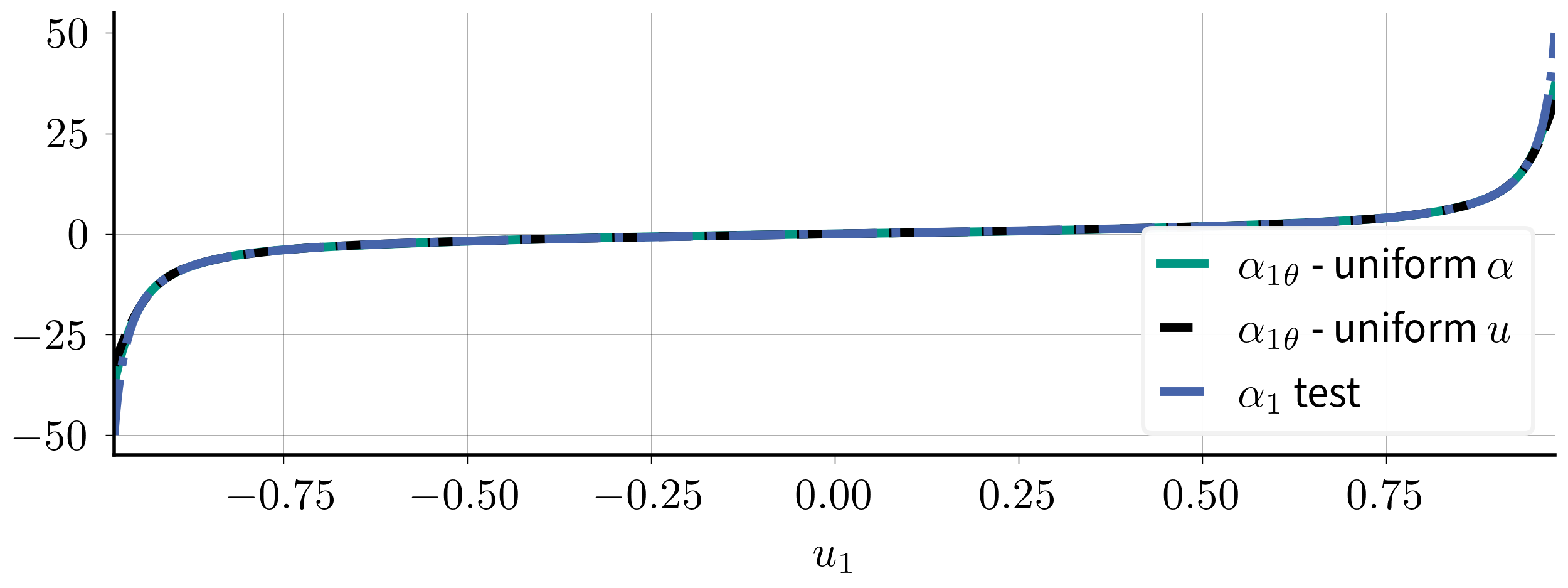}
            {\\a) $\alpha_{1\theta}$ and $\alpha_1$ over $u_1$}
    \end{minipage}
    \begin{minipage}{0.49\textwidth}
      \centering
       \includegraphics[width=\textwidth]{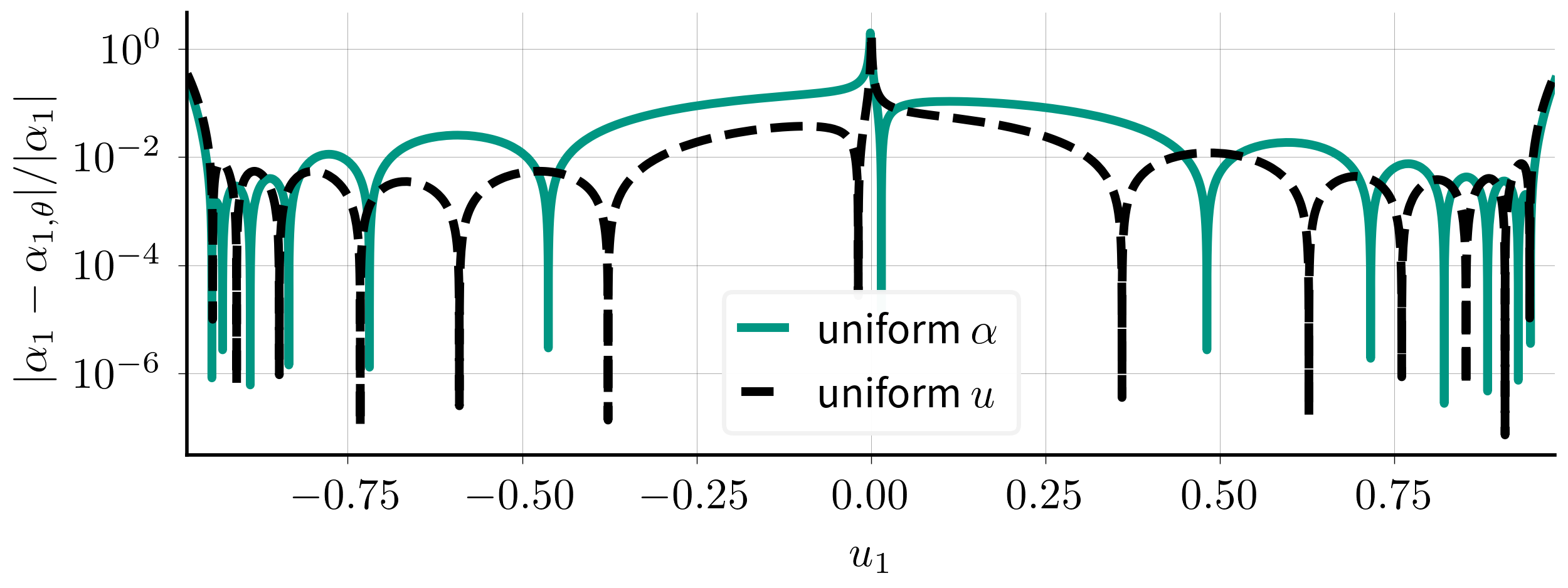}
            {\\b) Log difference in $\alpha_{1\theta}$ }
    \end{minipage}
     \begin{minipage}{0.49\textwidth}
       \centering
        \includegraphics[width=\textwidth]{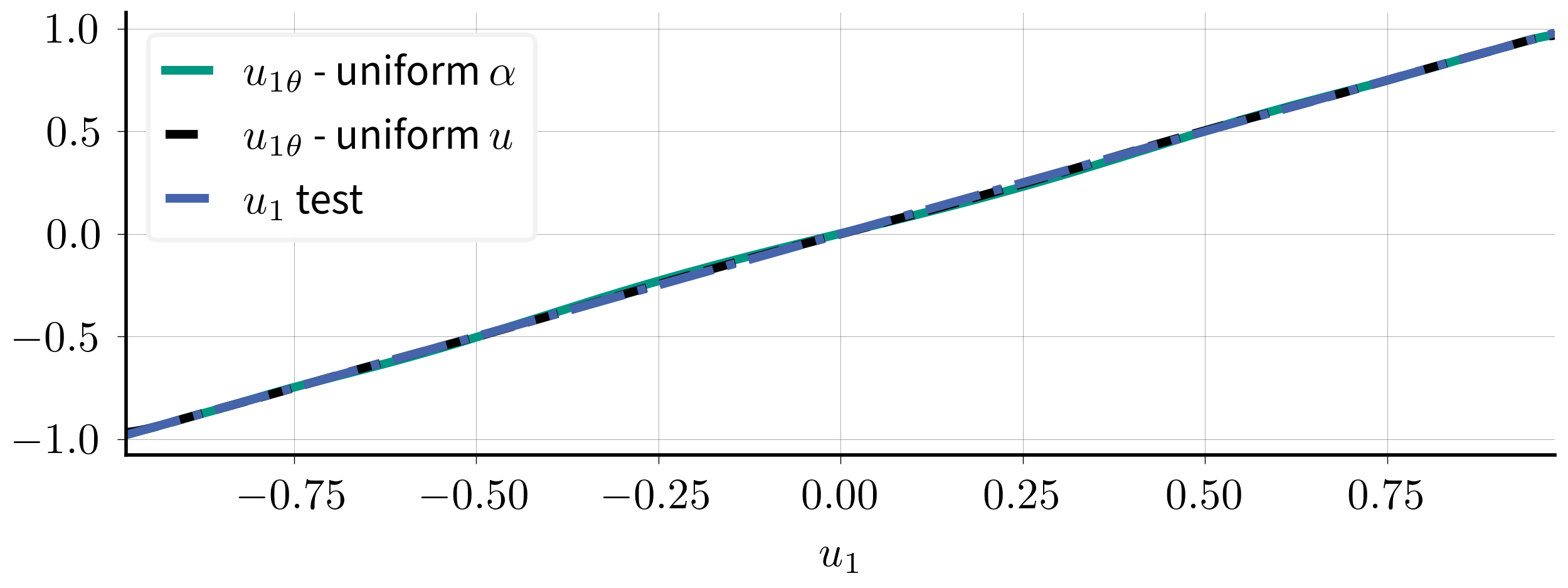}
            {\\c) Reconstructed $u_{1\theta}$ and $u_1$ over $u_1$}
    \end{minipage}
    \begin{minipage}{0.49\textwidth}
      \centering
       \includegraphics[width=\textwidth]{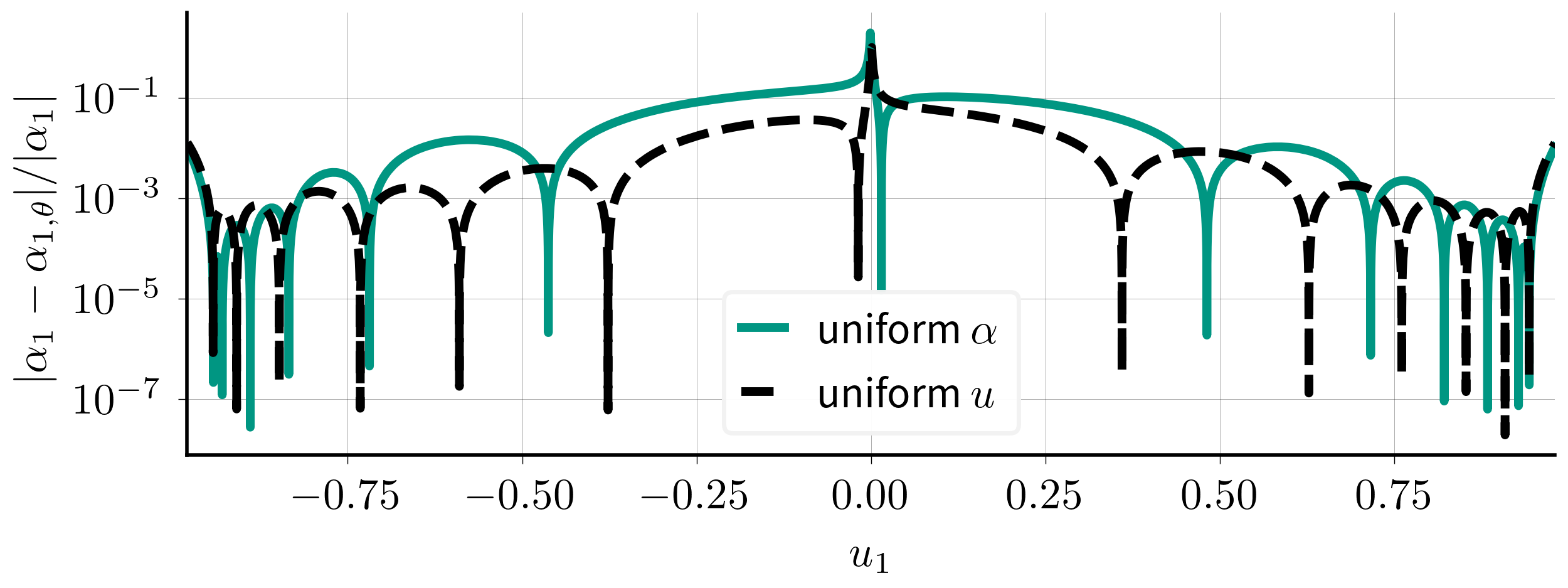}
            {\\d) Log difference in reconstructed $u_{1\theta}$}
    \end{minipage}
     \begin{minipage}{0.49\textwidth}
       \centering
        \includegraphics[width=\textwidth]{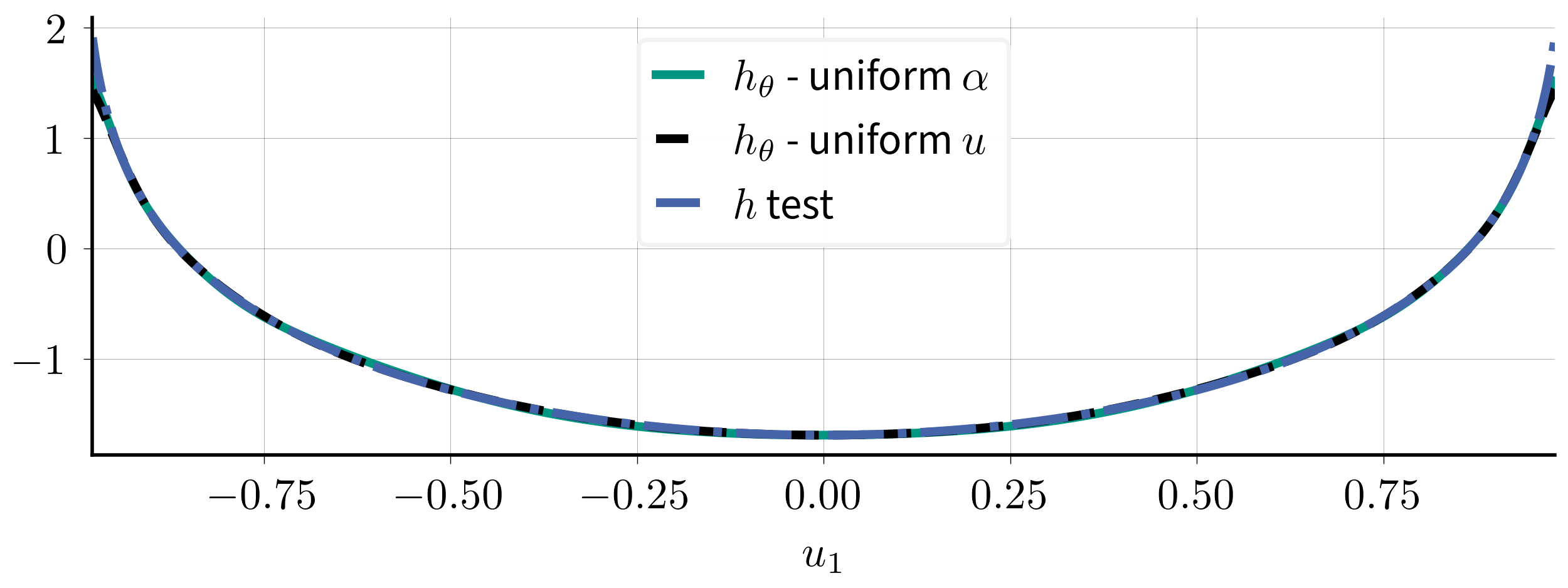}
            {\\e) $h_\theta$ and $h$ over $u_1$} 
    \end{minipage}
    \begin{minipage}{0.49\textwidth}
      \centering
       \includegraphics[width=\textwidth]{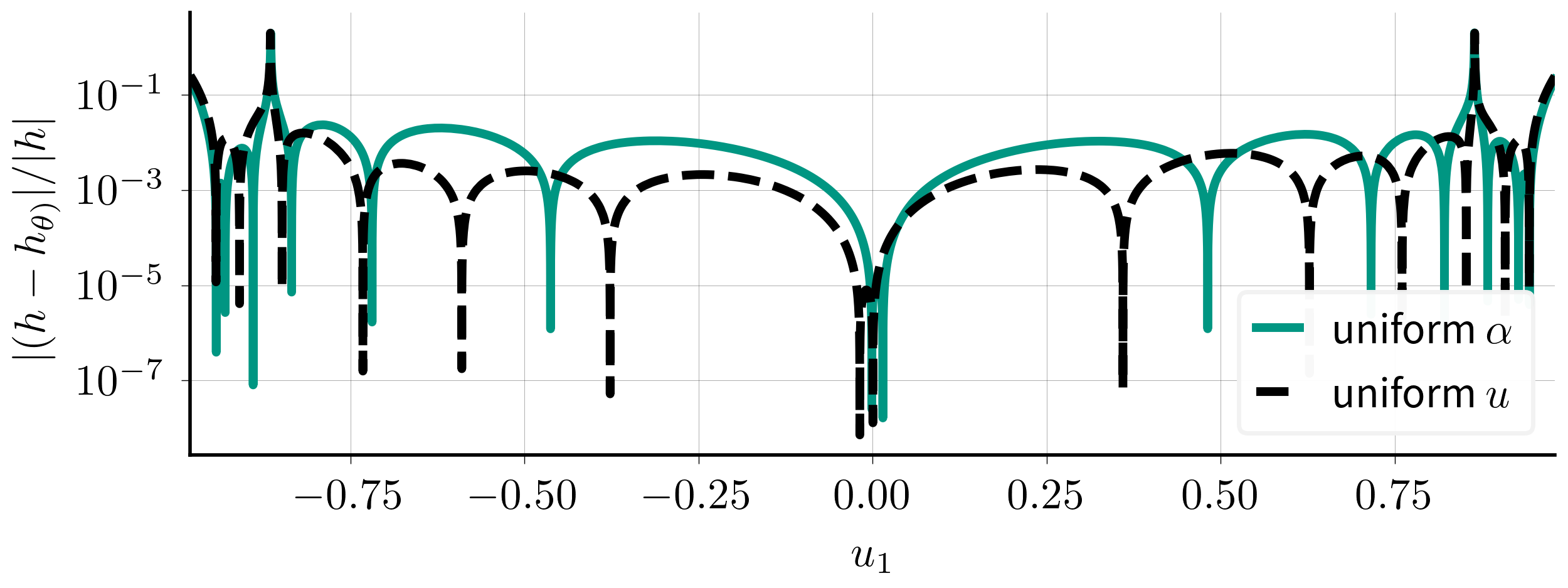}
            {\\a)  Log difference in $h_{\theta}$}
    \end{minipage}
    \caption{Validation tests for M$1$ $1$D networks trained on data-sets generated uniformly in $\alpha$ and $u$. $10^3$ training points and $10^4$ testing points. Distance to $\partial\mathcal{R}$ is $0.01$ }
    \label{fig_dataSampling1D}
\end{figure}
In this section, we compare networks trained on data using Algorithm~\ref{alg_samplingUniform} and Algorithm~\ref{alg_samplingUniformAlpha}. By construction, the data distribution generated by these strategies are completely different, as it is shown in Fig.~\ref{fig_dataSampling}. Figure~\ref{fig_dataSampling1D}a) shows, that as we approach the boundary $\partial\overline{\mathcal{R}}^r$, the values of $\alpha$ increase rapidly, while the slope is comparatively flat in the middle of the domain $[-1,1]$. Thus, uniform sampling in $\alpha$ results in a sparse distribution in the interior of $\overline{\mathcal{R}}^r$ and an accumulation of data-points near $\partial\overline{\mathcal{R}}^r$. In addition, it is not trivial to select the Lagrange multipliers in a way such that the distance of the set of corresponding moments is uniform to $\partial\overline{\mathcal{R}}^r$. This can be seen in the left image Fig.~\ref{fig_dataSampling}, where we sample $\alpha$ from a set $A=[-50,50]\times[-50,50]$. The upper boundary of the set of corresponding $\overline{u}^r$ is then concave, however the upper boundary of $\overline{\mathcal{R}}^r$ is a straight line from $(\overline{u}^r_1,\overline{u}^r_2)=(-1,1)$ to $(1,1)$.

We perform a synthetic validation test of identical models for the $M_1$ $1$D closure once trained on data generated with Algorithm~\ref{alg_samplingUniform} and once trained on data generated with Algorithm~\ref{alg_samplingUniformAlpha}. Since in this test case, $\mathcal{R}$ is one dimensional and we can sample the exact same interval with both strategies.  We train $50$ models for each data set using $10000$ training samples until convergence and compare the models. Figure~\ref{fig_dataSampling1D} shows the evaluations of the model using $100000$ unseen test data-points uniformly distributed in $\mathcal{R}$. The left column shows the evaluations $u_\theta$, $\alpha_\theta$, $h_\theta$ of both models and the test data $u$, $\alpha$ and $h$. The right column shows relative errors of the networks predictions. As expected, Algorithm~\ref{alg_samplingUniform} performs better in the interior of $\mathcal{R}$ in objective functions $h$, $\alpha$ and $u$, since in this region, there is more training data available. However, the advantage of Algorithm~\ref{alg_samplingUniformAlpha} near the boundary is not as big as the advantage of Algorithm~\ref{alg_samplingUniform} in the interior, where we experience almost an order of accuracy difference.
Note, that the error in $h_\theta$ decreases, as $u_1\rightarrow 0$, whereas the error in $u_\theta$ and $\alpha_\theta$ increases in this region. By Eq.~\eqref{eq_alpha0_recons} and the definition of $u_\theta$, the error in $u_\theta$ is nonlinearly dependent on the error in $\alpha_\theta$. The spike in the relative error is due to the proximity of $\alpha$ (and) $u$ to $0$. Similar behavior can be seen at the roots of $h$. 
As a conclusion, we use Algorithm~\ref{alg_samplingUniform} to generate training data for the models employed in the following test cases. 

\begin{table}[htbp]
\centering
  	\caption{Validation losses for different closure networks} 
	\begin{tabular}{|l|l|l|l|l|l|l|l|} 
		\hline
		 & Layout &  MSE$(h,h_\theta)$ & MSE$(\alpha,\alpha_\theta)$ & MSE$(u,u_\theta)$&   MAE$(h,h_\theta)$ & MAE$(\alpha,\alpha_\theta)$ & MAE$(u,u_\theta)$\\
		\hline 
		$M_1$ $2$D  & $18\times 8$ & $1.10\mathrm{e}{-6}$ & $3.39\mathrm{e}{-5}$ & $2.93\mathrm{e}{-6}$ & $1.02\mathrm{e}{-3}$ & $3.36\mathrm{e}{-3}$ & $0.96\mathrm{e}{-3}$ \\ 
		\hline 
		$M_1$ $1$D  & $10\times 7$ & $7.87\mathrm{e}{-7}$ & $7.52\mathrm{e}{-4}$ & $1.47\mathrm{e}{-6}$ & $7.57\mathrm{e}{-4}$ & $1.26\mathrm{e}{-2}$ & $9.59\mathrm{e}{-4}$ \\ 
		\hline
		$M_2$ $1$D  & $15\times 7$  & $1.33\mathrm{e}{-5}$ & $2.81\mathrm{e}{-4}$ & $2.81\mathrm{e}{-4}$ & $3.11\mathrm{e}{-3}$ & $1.23\mathrm{e}{-2}$ & $1.23\mathrm{e}{-2}$ \\  	\hline
	\end{tabular}  \label{tab_val_losses}
\end{table}
All networks are build using the input convex architecture discussed in Section~\ref{sec_neuralEntropyClosure}. Each network has a dense input layer followed by a block of convex bridge layers defined in Eq.~\eqref{eq_ICnnKons}. This convex block is described by the layout column in Table~\ref{tab_val_losses}, with format width $\times$ depth. After the convex block, a convex bridge layer with half the width of the block followed by a convex output layer is added. The output layer has the same design as a convex bridge layer, but lacks the softplus activation function, since we deal with a regression task. 
The validation loss of the neural networks after training can be seen in Table~\ref{tab_val_losses}. The networks are trained on an Nvidia RTX 3090 GPU in single-precision floating-point accuracy. Notice, that the mean absolute error in $\alpha$ is significantly higher than the error in $h$ or $u$ for all networks. The reason for this is the high range of values, that $\alpha$ can attain, which is inconvenient for the neural network.

\subsection{Computational Efficiency}
In the following, we compare the computational efficiency of the neural network surrogate model and the Newton optimizer in an isolated, synthetic test case. We consider the $M_2$ closure in $1D$ and use only normalized moments. In contrast to the neural network, the performance of the Newton solver is dependent on the proximity of the moments $u$ to the boundary $\partial\overline{\mathcal{R}}$, we consider three test cases. First, the moments are uniformly sampled using Algorithm~\ref{alg_samplingUniform}, second we only draw moments near the center of $\overline{\mathcal{R}}$ and lastly, we use only moments in proximity to  $\partial\overline{\mathcal{R}}$. The Newton solver is implemented in the KiT-RT~\cite{KITRT} framework. In the kinetic scheme, there is one optimization problem for each grid cell. Furthermore, optimization problems of different grid cells are independent of each other. A meaningful and straight-forward way to parallelize the minimal entropy closure is to employ one instance of the Newton optimizer per available CPU-core that handles a batch of cells. For comparability, we set the accuracy tolerance of the Newton solver to single-precision floating point accuracy. On the other hand, we interpret the number of grid cells as the batch size in of the neural entropy closure. Parallelization is carried out by the tensorflow backend.

We execute the neural network once on the CPU and once on  the GPU. The used CPU is a $24$ thread AMD Ryzen9 3900x with 32GB RAM and the GPU is a RTX3090 with 20GB RAM. The experiments are reiterated $100$ times to reduce time measurement fluctuations. Table~\ref{tab_timingbenchmark} displays the mean timing for each configuration and corresponding standard deviation.
\begin{table}[htbp]
	\caption{Computational cost for one iteration of the 1D solver in seconds $s$} 
	\centering
	\begin{tabular}{|l|l|l|l|} 
		\hline
		 & Newton & neural closure CPU & neural closure GPU \\
		\hline 
		uniform, $10^3$ samples & $0.00648\pm 0.00117$ s & $0.00788\pm 0.00051$ s & $0.00988\pm 0.00476$ s \\ 
		\hline 
		uniform, $10^7$ samples & $5.01239\pm 0.01491$ s & $0.63321\pm 0.00891$ s  & $0.03909\pm 0.00382$ s \\
		\hline
		boundary, $10^3$ samples & $38.35292\pm  0.07901$ s & $0.00802\pm 0.00064$ s & $0.00974\pm 0.00475$ s\\
		\hline 
		boundary, $10^7$ samples &  $27179.51012\pm 133.393$  s & $0.63299\pm 0.00853$ s  & $0.03881\pm 0.00352$ s  \\  
		\hline
		interior, $10^3$ samples & $ 0.00514\pm 0.00121$ s & $0.00875\pm 0.00875$ s & $0.00956\pm  0.00486$ s\\
		\hline 
		interior, $10^7$ samples & $4.24611\pm 0.03862$ s & $0.63409\pm 0.00867$ s  & $0.03846 \pm 0.00357$ s\\
		\hline
	\end{tabular} 
	\label{tab_timingbenchmark}
\end{table}
Considering Table~\ref{tab_timingbenchmark}, we see that the timing of a neural network is independent on the condition of the optimization problem, whereas the Newton solver is $6300$ times slower on a on a moment $u$ with $\norm{u-\partial\overline{\mathcal{R}}}_2 =0.01$ compared to a moment in the interior. The average time to compute the Lagrange multiplier of a uniformly sampled moment $u$ is $27$\% higher than a moment of the interior.
However, we need to take into account that the neural entropy closure is less accurate near $\partial\overline{\mathcal{R}}$ as shown in Fig.~\ref{fig_dataSampling1D}.
Furthermore, we see that the acceleration gained by usage of the  neural network surrogate model is higher in cases with more sampling data. This is apparent in the uniform and interior sampling test cases, where the computational time increases by a factor of $\approx 73$, when the data size increases by a factor of $10^4$. The time consumption of the Newton solver increases by a factor of $\approx840$ in the interior sampling case, respectively $\approx 782$ in the uniform sampling case. Note, that in this experiment, all data points fit into the memory of the GPU, so it can more efficiently perform SIMD parallelization.


\begin{table}[htbp]
    \centering
	\caption{Computational setup of the test cases} 
	\begin{tabular}{|l|l|l|l|l|} 
		\hline
		&1D M1 Linesource & 1D M2 Linesource & 2D M1 Periodic \\
		\hline
		$T$ & $(0,0.7]$ & $(0,0.7]$ & $(0,1.84]$ \\	\hline
		$N_t$  & 1399 & 1399 & $1226$\\ 	\hline
		$X$ & $[0,1]$ & $[0,1]$& $[-1.5,1.5]\times[-1.5,1.5]$ \\	\hline
		$N_X$ &$100$ & $100$ & $100^2$ \\	\hline
		Quadrature &  Legendre & Legendre& Tensorized Gauss Legendre\\	\hline
		$\mathbf{V}$ &  $[-1,1]$ & $[-1,1]$ & $[-1,1]\times[0,2\pi)$ \\ 	\hline
		$N_\mathbf{V}$ & $28$ & $28$ & $200$\\	\hline
		Basis & Monomial &  Monomial & Monomial\\ 	\hline
		CFL  & $0.05$ & $0.05$ &  $0.1$ \\ 	\hline
		$\sigma$  & $1.0$ & $1.0$ & $0.0$\\ 	\hline
	\end{tabular} 
	\label{tab_2D_setup}
\end{table}


\subsection{Inflow into purely scattering homogeneous medium in 1D}
\begin{figure}
    \centering
    \begin{minipage}{0.495\textwidth}
        \includegraphics[width=\textwidth]{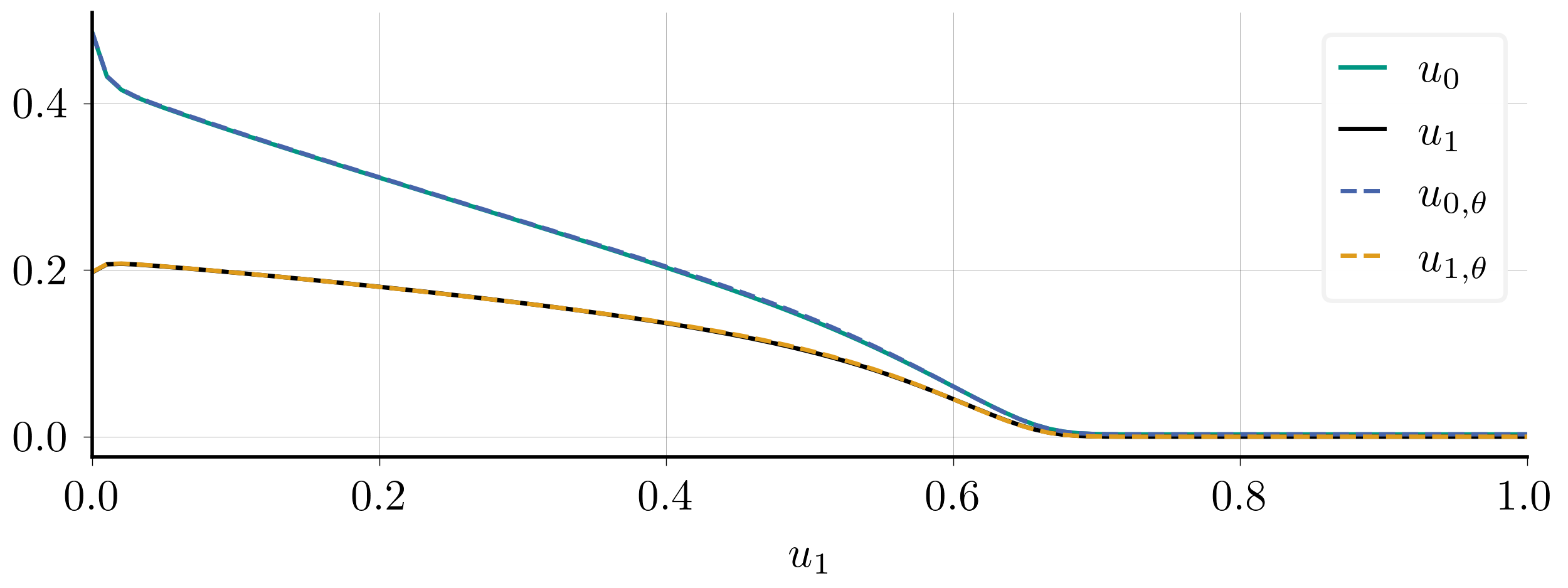}
   \end{minipage}
   \begin{minipage}{0.495\textwidth}
       \includegraphics[width=\textwidth]{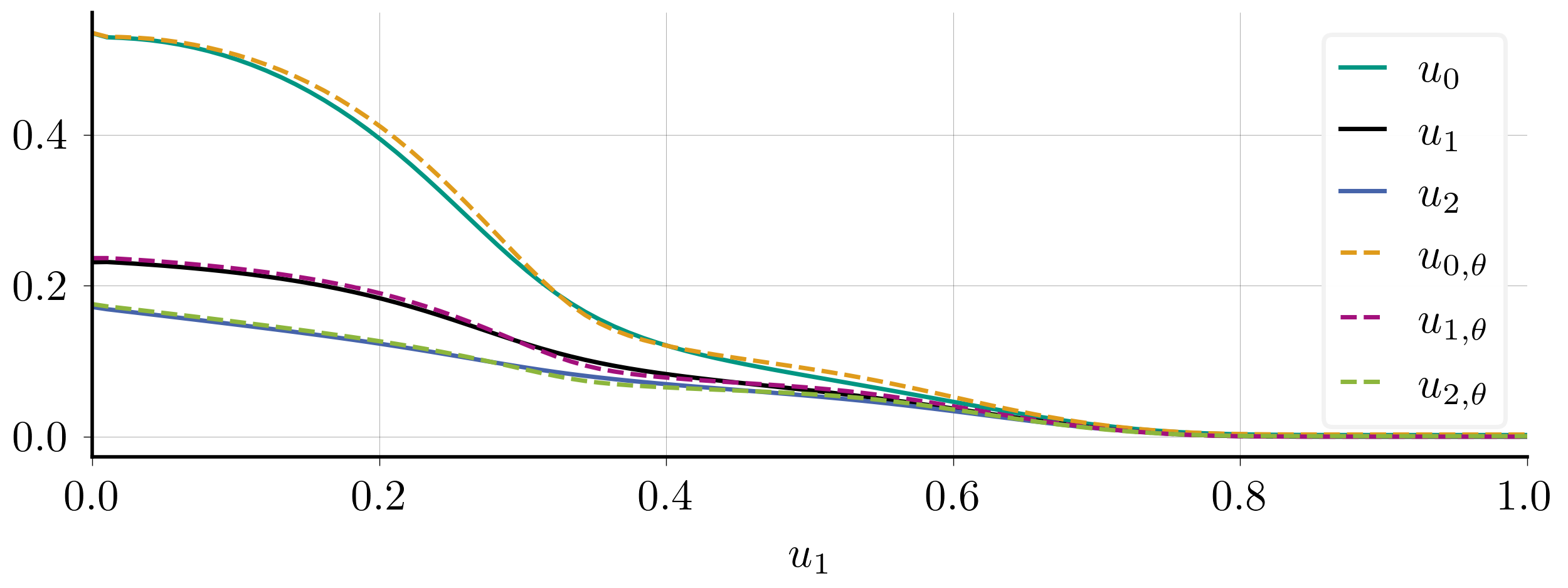}
   \end{minipage}
   \caption{Comparison of neural closed and benchmark solution at $t=0.7$. The left figure displays the $M_1$ and the right figure displays the $M_2$ test case.}
   \label{fig_1dm1}
\end{figure}

\begin{figure}
   \begin{minipage}{0.495\textwidth}
       \includegraphics[width=\textwidth]{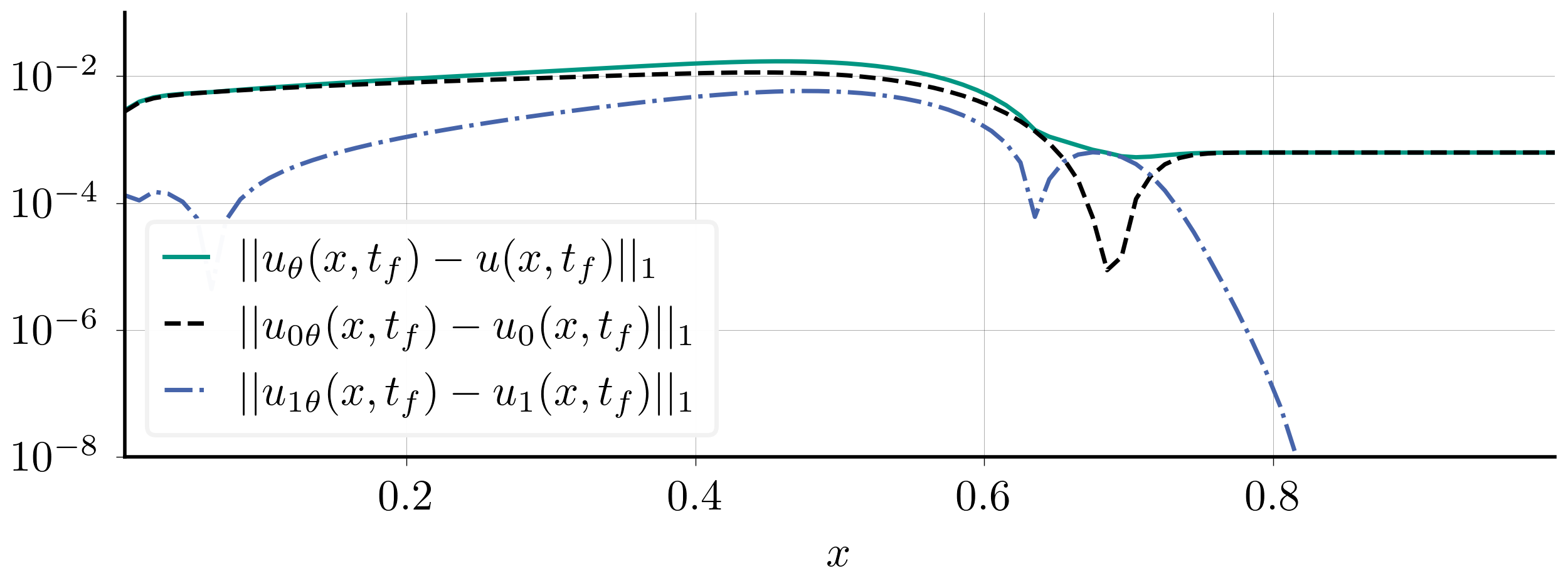}
    \end{minipage}
    \begin{minipage}{0.495\textwidth}
       \includegraphics[width=\textwidth]{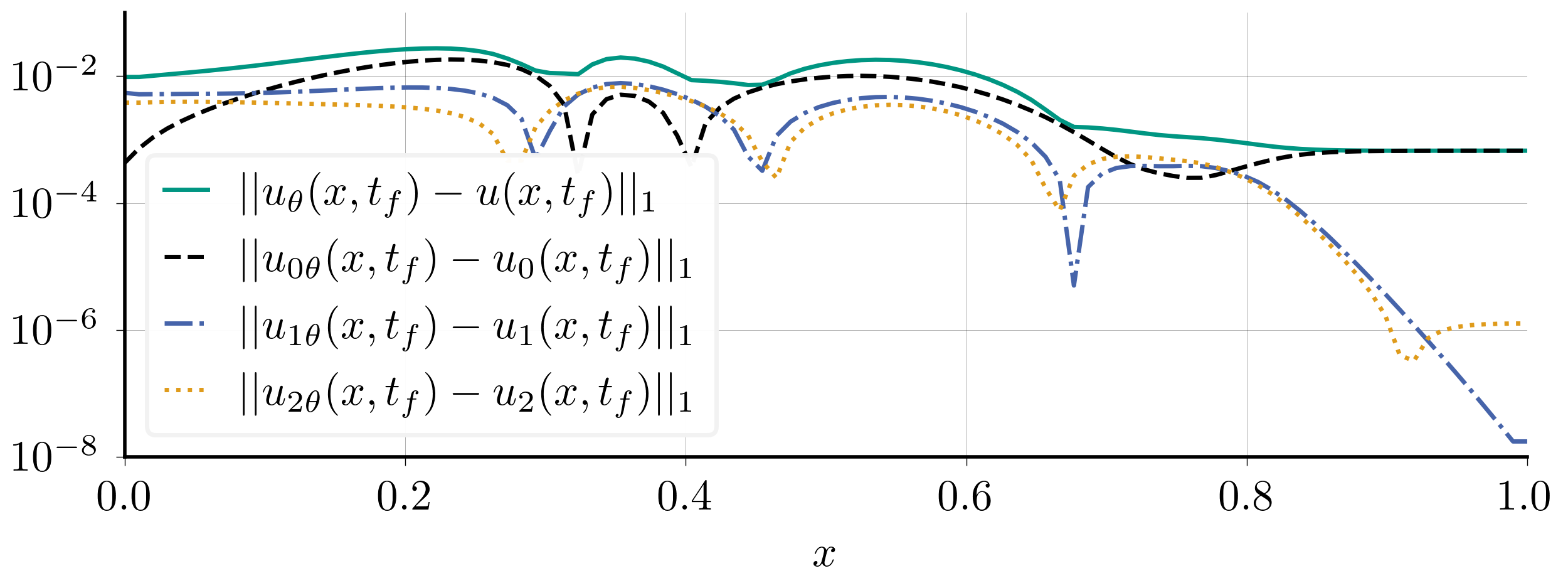}
    \end{minipage}
   \caption{Comparison of the systems entropy over time, $M_1$ case on the left and $M_2$ case on the right.}
   \label{fig_err}
\end{figure}

Let us first study the particle transports in an isotropic scattering medium.
We consider the one-dimensional geometry, where the linear Boltzmann equation reduces to
\begin{equation}
    \partial_t f+\mu \partial_{\mathbf{x}} f=Q(f)=\sigma\left(\frac{1}{2}\inner{ f} - f\right),
\end{equation}
and the corresponding moment model becomes
\begin{align}
  \partial_t u + \partial_{\mathbf{x}} \inner{\mu m f_u}&= \inner{m Q(f_u)}\\
      f_u &= \eta_*(\partial_u\mathcal{N}_\theta(u)\cdot m)
\end{align}
The initial domain is set as vacuum with $f(0,\mathbf{x},\mu)=0$,
and an inflow condition is imposed at the left boundary of domain with $f(t>0,0,\mu>0)=0.5$.
The detailed computational setup can be found in Table \ref{tab_2D_setup}.
The solution profiles at final time $t_f$ of the neural entropy closed moment system and the benchmark solver are presented in Fig.~\ref{fig_1dm1} for the $M_1$ and $M_2$ system. 
In Fig.~\ref{fig_err}, we see the corresponding errors of the $M_1$ and $M_2$ solution for each grid cell. The plots show the error of the moment vector $u$ in $l_1$ norm as well as the $l_1$ errors of its components. 
We notice a slight increase of the error in the $M_2$ case, overall it displays a good approximation as well.

\subsection{2D Test cases}
\begin{figure}
    \centering
   \begin{minipage}{0.49\textwidth}
        \includegraphics[width=\textwidth]{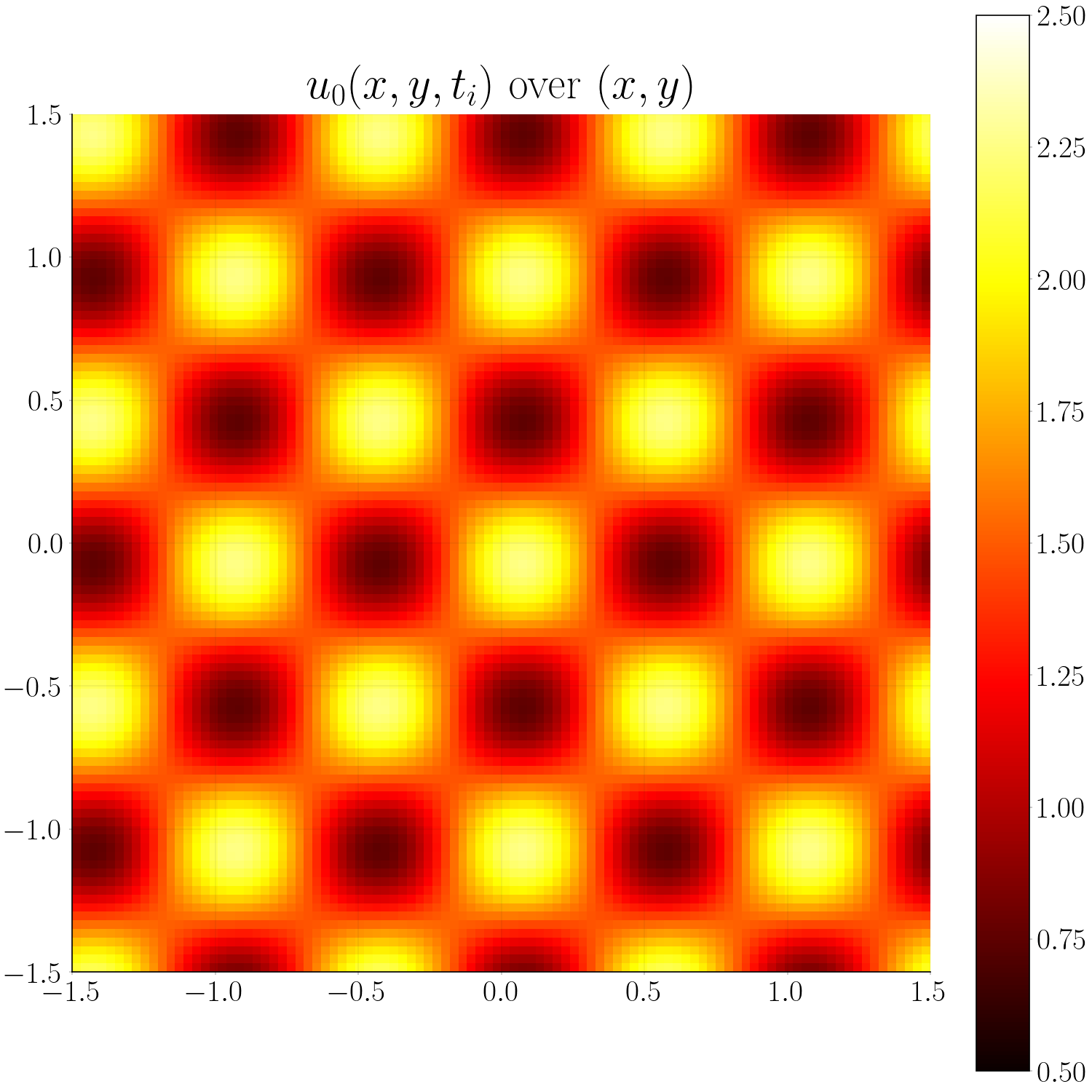}
   \end{minipage}
    \begin{minipage}{0.49\textwidth}
       \includegraphics[width=\textwidth]{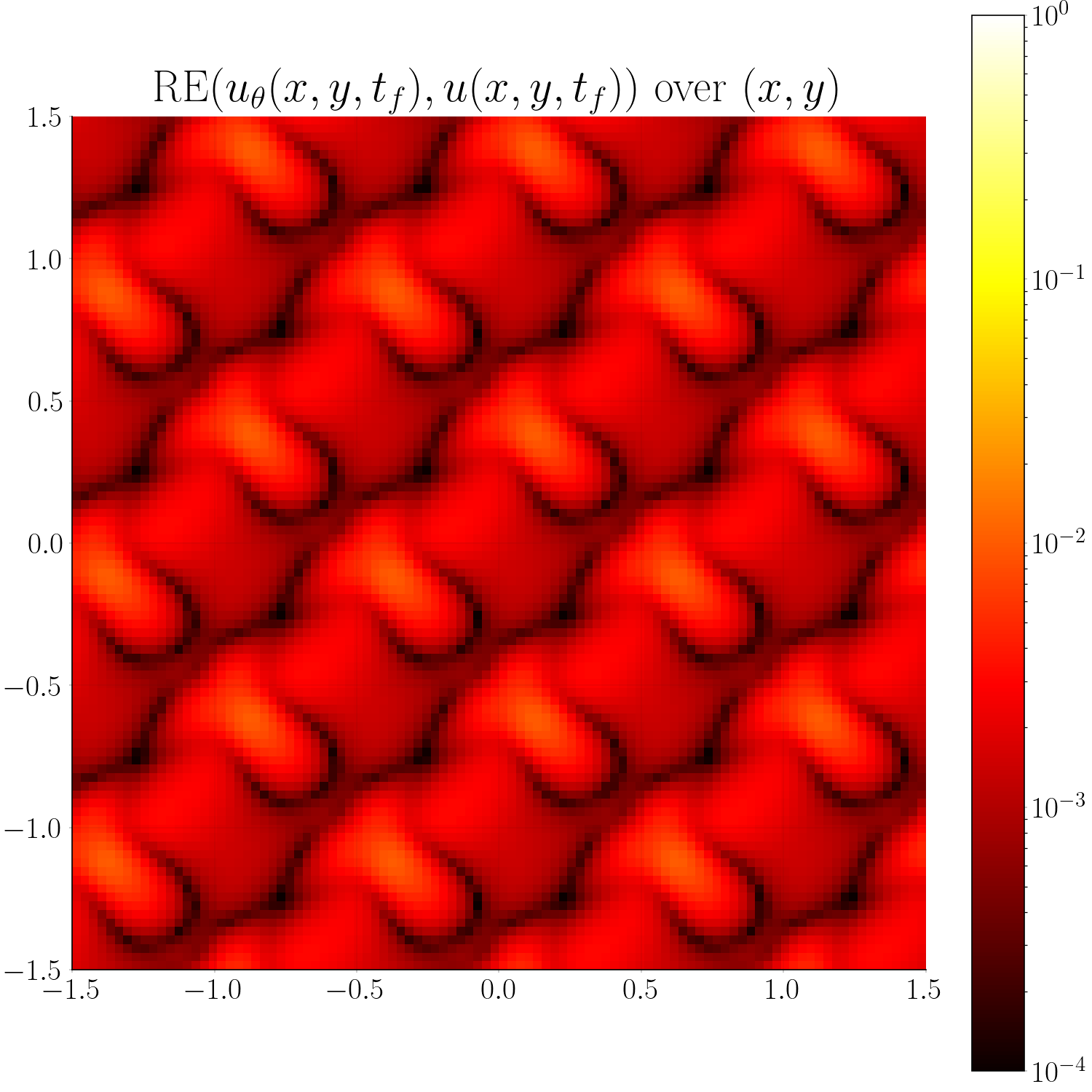}
   \end{minipage}
   \caption{Solution and relative errors in the periodic $2D$ test case at iteration $i=60$.}
  \label{fig_period2D}
\end{figure}

\begin{figure}
    \centering
      \begin{minipage}{0.49\textwidth}
            \includegraphics[width=\textwidth]{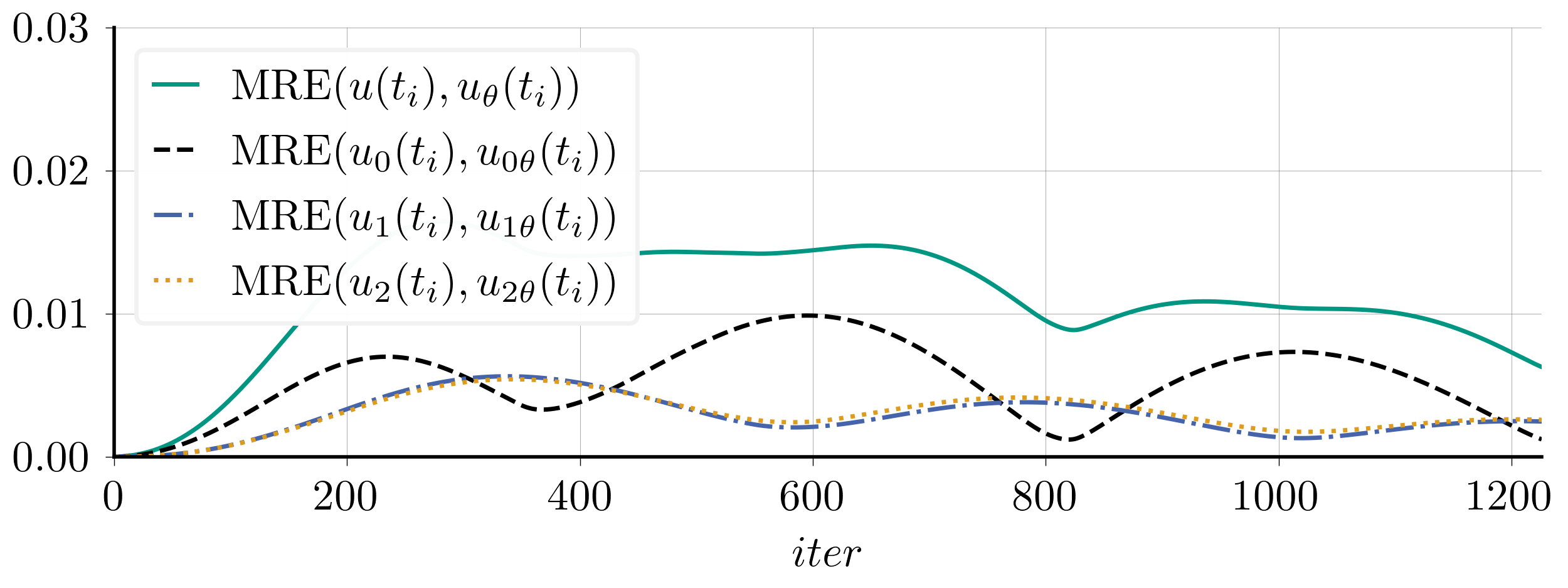}
   \end{minipage}
    \begin{minipage}{0.49\textwidth}
            \includegraphics[width=\textwidth]{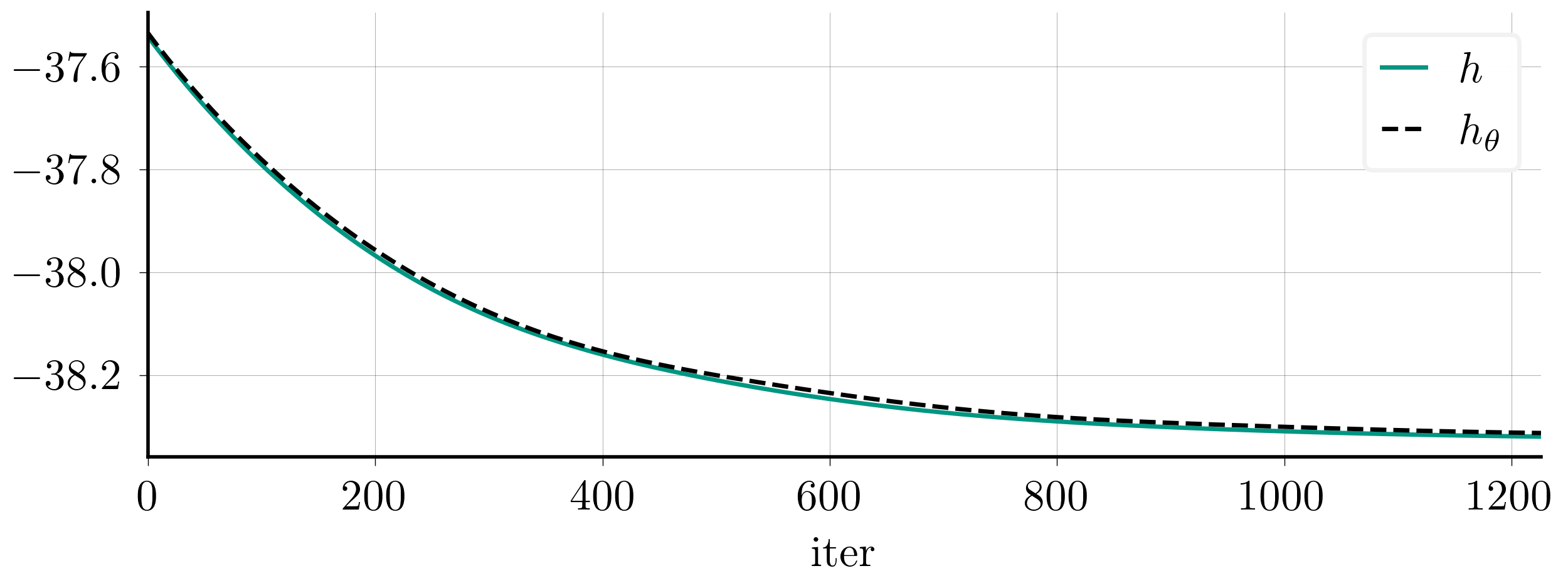}
    \end{minipage}
    \caption{ Mean relative error of neural entropy closed and Newton closed solution (left) system entropy (right) at each time step.}\label{fig_MRE_entropy}
\end{figure}

We consider a rectangular domain in two spatial dimensions. The phase space of the Boltzmann equation is thus five dimensional, where $\mathbf{X}=[-1.5,1.5]^2$, $\mathbf{V}=\menge{v\in\mathbb{R}^2:\norm{v}_2 < 1}$ and $t>0$. We consider the $M_1$ closure with a monomial basis $m(v)=[1,v_x,v_y]^T$, where $v_x$ and $v_y$ are defined in Table~\ref{tab_quad}.
The Boltzmann equation reduces to 
\begin{align}
\begin{aligned}
    \partial_t f + v_x\partial_x f + v_y\partial_y f =& Q(f)(v) \\
    =&\sigma\left(\inner{\frac{1}{2\pi}f} -f\right),
\end{aligned}
\end{align}
where $\sigma$ is the scattering coefficient. The corresponding moment system with neural entropy closure reads as
\begin{align}
\begin{aligned}
  \partial_t u +\partial_x\inner{v_xmf_u}++\partial_y\inner{v_y mf_u} &= \inner{mQ(f)} \\
    f_u &= \eta_*(\partial_u\mathcal{N}_\theta(u)\cdot m)
\end{aligned}
\end{align}
We inspect the $M_1$ closure choose the systems initial condition as
\begin{align}
    u_0 &= 1.5 + \cos(2\pi x)\cos(2\pi y), \qquad (x,y)\in\mathbf{X}\\
    u_1 &= 0.3 u_0 \\
    u_2 &= 0.3 u_0
\end{align}
The moment system is discretized using the kinetic scheme described in Section~\ref{sec_kineticSolver}.  We compare the kinetic scheme once closed using a Newton solver and once closed using the neural entropy closure. 
 We use the solver configuration in Table~\ref{tab_2D_setup} and an entropy consistent numerical scheme. We run the simulation until a final time $t_f=3.5$, which translates to $1630$ time-steps. 
Figure~\ref{fig_period2D} shows the neural entropy closed solution at time iteration $i=60$ on the left image and the relative error of the entropy closed solution compared to the Newton optimizer closed solution on the right image. The solution snapshot shows no signs of numerical artifacts and indeed, the relative error per grid cell is mostly in the order of $10^{-2}$ or lower.
In the second experiment, we run both solvers up to final time $t_f$ and compare the propagated error of the neural entropy closure to the Newton solution, that is configured to solve the entropy optimization problem Eq.~\eqref{eq_entropyDualOCP} up to machine precision. The results can be seen in Fig.~\ref{fig_MRE_entropy}. Note, that the mean relative error does not increase over a threshold of $2.3$\%. 

Lastly, we look at the entropy of the system at each time step. Due to the periodic boundaries, the physical system is closed. We have chosen the upwind scheme for the numerical flux of the moment system, which is an entropy dissipating scheme. Figure~\ref{fig_MRE_entropy} shows that the Newton and neural entropy closed system both dissipate mathematical entropy over time and furthermore, the dissipation rate of both schemes is similar. This shows, that the entropy dissipation property is conserved by the neural entropy closure.

\section{Summary and Conclusion}

In this paper we addressed the moment system of the linear Boltzmann equation, its minimal entropy closure, and the challenges of classical numerical approaches.
We introduced a novel convex neural network based approach to close the moment hierarchy of the linear Boltzmann equation. 
The nature of the entropy minimization problem allows clear definition of the convex set of all possible input data for the neural network. On the other hand, the problem is ill conditioned on the boundary of the realizable set and thus poses significant challenges for generating training data on and near the boundary of it. We propose two algorithms to tackle this challenge and compared them in terms of training and validation performance of the corresponding networks. We found that uniform sampling of $\mathcal{R}$ is advantageous to uniform sampling of the Lagrange multipliers $\alpha$. 

Next, we have constructed an input convex neural network, that mimics the entropy functional of the moment system of the Boltzmann equation. We successfully employed this neural entropy closure in several $1$D and $2$D test cases and systems of different moment order.
We found a good agreement between the neural entropy closed system and the solutions computed with a conventional Newton solver within the boundaries of the training performance of the neural networks. As expected, the neural network based closure is significantly more efficient in terms of computational time  compared to a Newton solver.

Several challenges remain to be studied. In higher spatial dimensions and higher order moment systems, it is not clear how to characterize $\overline{\mathcal{R}}$ directly, thus one needs to apply Algorithm~\ref{alg_samplingUniformAlpha} and sample in the space of Lagrange multipliers. Then, one needs to explore other sampling distributions due to their shown disadvantages compared to uniform sampling in $\overline{\mathcal{R}}$. Furthermore, it is preferable to characterize the distance to $\partial\overline{\mathcal{R}}$ in terms of $\alpha$ in order to construct an efficient sampling algorithm.
Additionally, a rigorous error analysis for the neural entropy closure needs to be conducted. 
Further research will be dedicated to a sophisticated sampling algorithm for higher dimensions, and consider an error analysis of the neural network reconstruction. 

\section*{Acknowledgements}
The authors acknowledge support by the state of Baden-Württemberg through bwHPC. Furthermore, the authors would like to thank Dr. Jonas Kusch for fruitful discussions about realizability and the minimal entropy closures as well as Jannick Wolters for support in scientific computing matters.
The research is funded by the Alexander von Humboldt Foundation (Ref3.5-CHN-1210132-HFST-P).

The work of Cory Hauck is sponsored by the Office of Advanced Scientific Computing Research, U.S. Department of Energy, and performed at the Oak Ridge National Laboratory, which is managed by UT-Battelle, LLC under Contract No. De-AC05-00OR22725 with the U.S. Department of Energy. The United States Government retains and the publisher, by accepting the article for publication, acknowledges that the United States Government retains a non-exclusive, paid-up, irrevocable, world-wide license to publish or reproduce the published form of this manuscript, or allow others to do so, for United States Government purposes. The Department of Energy will provide public access to these results of federally sponsored research in accordance with the DOE Public Access Plan (http://energy.gov/downloads/doe-public-access-plan).
\bibliography{sample}

\end{document}